\documentclass[1p]{elsarticle}

\usepackage{lineno,hyperref}
\modulolinenumbers[5]

\journal{Journal of Geometry and Physics}

\usepackage{graphicx}
\usepackage[inline]{enumitem}
\usepackage{subcaption} 
\usepackage{colortbl}
\usepackage{float}

\usepackage{natbib}
\usepackage{bibentry}

\usepackage{multirow}
\usepackage{hyphenat}
\usepackage[ruled, vlined, linesnumbered]{algorithm2e}
\usepackage{rotating} 
\usepackage{adjustbox}%
\usepackage[11pt]{moresize}
\usepackage{gensymb} 

\usepackage{amsthm}
\usepackage{amsmath}
\usepackage{amssymb}
\usepackage{amsfonts}

\usepackage{stmaryrd} 

\usepackage[normalem]{ulem}

\usepackage{csquotes}

\usepackage{mathrsfs}

\newtheorem{theorem}{Theorem}[section]

\newtheorem{proposition}[theorem]{Proposition}

\theoremstyle{definition}

\newtheorem{remark}[theorem]{Remark}

\hypersetup{hidelinks}

\begin{document}

\begin{frontmatter}

\title{Port-Hamiltonian Modeling of Ideal Fluid Flow: \\Part I. Foundations and Kinetic Energy}

\author{Ramy Rashad\textsuperscript{1,*}\corref{RAM}, Federico Califano\textsuperscript{1}, Frederic P. Schuller\textsuperscript{2}, Stefano Stramigioli\textsuperscript{1}}
\address{\textsuperscript{1} Robotics and Mechatronics Department, University of Twente, The Netherlands}
\address{\textsuperscript{2} Department of Applied Mathematics, University of Twente, The Netherlands}
\cortext[RAM]{Corresponding author. Email: r.a.m.rashadhashem@utwente.nl}

%
%

\begin{abstract}
In this two-parts paper, we present a systematic procedure to extend the known Hamiltonian model of ideal inviscid fluid flow on Riemannian manifolds in terms of Lie-Poisson structures to a port-Hamiltonian model in terms of Stokes-Dirac structures.
The first novelty of the presented model is the inclusion of non-zero energy exchange through, and within, the spatial boundaries of the domain containing the fluid.
The second novelty is that the port-Hamiltonian model is constructed as the interconnection of a small set of building blocks of open energetic subsystems.
Depending only on the choice of subsystems one composes and their energy-aware interconnection, the geometric description of a wide range of fluid dynamical systems can be achieved.
The constructed port-Hamiltonian models include a number of inviscid fluid dynamical systems with variable boundary conditions. Namely, compressible isentropic flow, compressible adiabatic flow, and incompressible flow.
Furthermore, all the derived fluid flow models are valid covariantly and globally on n-dimensional Riemannian manifolds using differential geometric tools of exterior calculus.
\end{abstract}


\end{frontmatter}


\newcommand{\half}{\frac{1}{2}}

\newcommand{\B}[1]{\boldsymbol{#1}}
\newcommand{\bb}[1]{\mathbb{#1}}
\newcommand{\cl}[1]{\mathcal{#1}}
\newcommand{\Rn}{\bb{R}^n}

\newcommand{\TwoTwoMat}[4]{
\begin{pmatrix}
#1 & #2 \\
#3 & #4
\end{pmatrix}}
\newcommand{\TwoVec}[2]{
\begin{pmatrix}
#1\\
#2 
\end{pmatrix}}

\newcommand{\ThrVec}[3]{
\begin{pmatrix}
#1\\
#2\\
#3
\end{pmatrix}}

\newcommand{\map}[3]{{#1}:{#2}\rightarrow {#3}}
\newcommand{\fullmap}[5]{
\begin{split}
#1 : {#2} &\rightarrow {#3}\\
	{#4} &\mapsto {#5}
\end{split}
}
\newcommand{\pair}[2]{\left\langle \left.  #1  \right|  #2 \right\rangle}

\newcommand{\blue}{\color{blue}}

\newcommand{\pH}{port-Hamiltonian }

\newcommand{\gothg}{\mathfrak{g}}
\newcommand{\gothgV}{\mathfrak{g^v}}
\newcommand{\gothgstar}{\mathfrak{g}^{*}}

\newcommand{\extd}{\textrm{d}}

\newcommand{\goths}{\mathfrak{s}}
\newcommand{\gothsStar}{\mathfrak{s}^{*}}

\newcommand{\secTanBdl}[1]{\Gamma(T#1)}

\newcommand{\volF}{\mu_{\normalfont\text{vol}}}

\newcommand{\abar}{{\bar{a}}}


\newcommand{\LieD}[2]{\cl{L}_{#1}#2}
\newcommand{\dt}{\frac{d}{dt}}

\newcommand{\dtLine}[1]{\frac{d}{dt}\bigg|_{#1}}
\newcommand{\depsLine}[1]{\frac{d}{d \epsilon}\bigg|_{#1}}

\newcommand{\pt}{\frac{\partial}{\partial t}}
\newcommand{\JLBrack}[3]{\llbracket{ #2},{#3} \rrbracket_{#1}}

\newcommand{\inner}[2]{\left \langle #1 , #2 \right\rangle}

\newcommand{\Lbrack}[2]{\left[#1,#2 \right]}
\newcommand{\LbrG}[2]{\Lbrack{#1}{#2}_{\gothg}}
\newcommand{\LbrS}[2]{\Lbrack{#1}{#2}_{\goths}}

\newcommand{\Pbrack}[2]{\left \{#1,#2 \right\}}

\newcommand{\parXi}[1][i]{\frac{\partial}{\partial x^{#1}}}
\newcommand{\dXi}[1][i]{\extd{ x^{#1}}}

\newcommand{\varD}[2]{\frac{\delta {#1}}{\delta {#2}}}

\newcommand{\adS}[1]{\boldsymbol{ad}_{#1}}
\newcommand{\adSdual}[1]{\boldsymbol{ad}^*_{#1}}

\newcommand{\spKForm}[2]{\Omega^{#1}(#2)}
\newcommand{\spVecF}[1]{\Gamma(T#1)}
\newcommand{\spVecX}[2]{\mathfrak{X}_{#2}(#1)}
\newcommand{\spKFormM}[1]{\spKForm{#1}{M}}
\newcommand{\spKFormMbc}[2]{\Omega^{#1}_{#2}(M)}
\newcommand{\spVecM}{\spVecX{M}{}}
\newcommand{\spVecMt}{\spVecX{M}{T}}
\newcommand{\spFn}[1]{C^\infty(#1)}

\newcommand{\vMeas}[1]{\textbf{v}(#1)}
\newcommand{\mMeas}[1]{\textbf{m}(#1)}
\newcommand{\mForm}{\mu_t}
\newcommand{\divr}[1]{\text{\normalfont div}(#1)}
\newcommand{\sdiff}{D_\mu(M)}
\newcommand{\RVF}[1]{\cl{R}(#1)}
\newcommand{\connect}[1]{\stackrel{\bb{#1}}{\nabla}}
\newcommand{\adIflat}[2]{ad_{#1}^*(\bb{I}^\flat {#2})}
\newcommand{\vf}{{\tilde{v}}}
\newcommand{\uf}{\tilde{u}}
\newcommand{\omV}{{\hat{\omega}}}

\newcommand{\mNmOne}{(-1)^{n-1} }

\newcommand{\effOne}{\delta_{\mu}H}
\newcommand{\effTwo}{\delta_{\alpha}H}
\newcommand{\bound}[1]{{#1}\rvert_{ \partial M}}
\newcommand{\effOneB}{\bound{\delta_{\mu}H}}
\newcommand{\effTwoB}{\bound{\delta_{\alpha}H}}

\newcommand{\gtInv}{g^{-1}_t}
\newcommand{\diffG}{\mathcal{D}(M)}

\newcommand{\parL}[1]{\frac{\delta l}{\delta {#1}}}

\newcommand{\Mbound}{\partial M}

\newcommand{\Ltwo}{\cl{L}_2}
\newcommand{\Mflat}{\bb{M}^\flat}
\newcommand{\Msharp}{\bb{M}^\sharp}

\newcommand{\state}{(\alpha,\mu)}

\tableofcontents

\section{Introduction}
Fluid mechanics is one of the most fundamental fields that has stimulated many ideas and concepts that are central to modern mathematical sciences.
Fluid mechanics has been studied in the literature using both the Lagrangian and Hamiltonian formalism.
In the classical Hamiltonian theory for fluid dynamical systems, a fundamental difficulty arises in incorporating the spatial boundary conditions of the system, which is also the case for general distributed parameter systems.
Previous Hamiltonian formulations of fluid flow in the literature \citep{Marsden1970,marsden1984semidirect,marsden1984reduction,Morrison1998} tend to focus on conservative systems with no energy-exchange with its surrounding environment.
Usually, this is imposed by certain assumptions on the system variables.
For example, if the spatial domain is non-compact, it is assumed that the system variables decay at infinity. Whereas if the spatial domain is compact, the boundary is assumed impermeable by imposing that the velocity vector field is tangent to the boundary.

Consequently, the traditional Hamiltonian theory is limited to distributed parameter systems on spatial manifolds without a boundary or ones with zero-energy exchange through the boundary.
While this is useful for analyzing a system that is isolated from its surroundings, it is certainly an obstacle for practical applications such as simulation and control.

The \pH framework grew out of the interest to extend the applicability of traditional Hamiltonian theory and has proven to be very effective for the treatment of both lumped-parameter \citep{maschke1992intrinsic} and distributed-parameter systems \citep{van2002hamiltonian}, see also the recent survey \citep{rashad2020twenty}.
While the traditional Hamiltonian formalism, in its generalized version on Poisson manifolds, is limited to conservative closed systems, the \pH formalism, based on Dirac structures, is applicable to non-conservative open systems capable of energy exchange with its environment.

Another core feature that distinguishes the \pH framework from the Hamiltonian framework is that a dynamical system is modeled as the interconnection of several subsystems classified based on their relation to energy.
Namely, energy storage, (free) energy dissipation, energy supply, and energy routing elements.
Therefore, the \pH framework is \textit{not a trivial extension} of the Hamiltonian framework but rather a \textit{paradigm shift} in treating dynamical systems in a bottom-up approach as opposed to the top-down approach of classical Hamiltonian theory.

\begin{figure}
\centering
\includegraphics[width=0.75\textwidth]{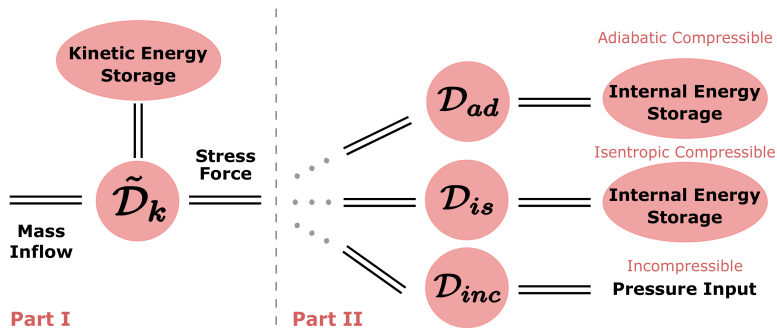}
\caption{Port-Hamiltonian model of various fluid dynamical systems as a network of interconnected energetic subsystems. Part I of this series deals with the kinetic energy storage $H_k$, its corresponding Dirac structure $\tilde{\cl{D}}_k$, and its open ports. Part II will dead with the various possibilities to interconnect another subsystem for storage of internal energy $H_i$ for adiabatic or isentropic compressible flow or incompressible flow.}
\label{fig:overview_picture}
\end{figure}

This two-part series of articles describes how fluid dynamical systems are systematically and completely modeled in the \pH framework by a small set of building blocks of open subsystems.
Depending on the choice of subsystems one composes in an energetically consistent manner, the geometric description of a number of fluid dynamical systems can be achieved, ranging from incompressible to compressible flows, as shown in Fig. \ref{fig:overview_picture}.

This decomposed network-based model of fluid dynamical systems comes at a great technical advantage.
Each of the subsystems is described in the structurally simplest possible way, even if other subsystems require a considerably more sophisticated formulation. 
The composition of such unequal subsystems is mediated by a Dirac structure, which routes the energy 
flow between all subsystems.

A second tremendous technical simplification is achieved concerning the precise choice of state space underlying each energetic subsystem.
For example when treating the kinetic energy of the fluid, no prior assumptions are required on the state space of fluid velocity to handle specific cases, such as impermeable boundaries and incompressiblity.
Instead, one chooses a state space that treats the general case of a compressible velocity field on any spatial domain (possible curved) that is permeable.
Then, a modeling assumption is imposed by composition with a suitable type of subsystem that physically models this assumption and is simply coupled to the unchanged kinetic-energy system.
The equations resulting from this procedure are, in the end, of course the same as those following from the traditional view, but their geometric description is technically simplified and made physically more insightful at the same time.

In addition to the system theoretic advantages of modeling ideal fluid flow in the \pH paradigm, the work we present here serves as a stepping stone for modeling fluid-structure interaction in the quest of understanding the flapping-flight of birds within the PORTWINGS project\footnote{http://www.portwings.eu/}.
In this project, we target a methodology which allows to describe the behavior of multi-physical systems, one of which is fluid dynamics, and the interconnection of various components in order to create a intrinsically composable model of the overall system.

In this first paper of the two-part series, we start the \pH modeling process of fluid dynamical systems with the construction of the energetic subsystem that stores the kinetic energy of a fluid flowing on a general spatial domain with permeable boundary.
The procedure to construct the \pH model we are aiming for relies greatly on understanding the underlying geometric structure of the state space of each energetic subsystem.
This geometric formulation, pioneered by \cite{arnold1966geometrie} and \cite{ebin1970groups}, will allow a systematic derivation of the underlying Hamiltonian dynamical equations and Dirac structures, usually postulated a priori in the literature \citep{van2001fluid,van2002hamiltonian,Polner2014,Altmann2017}.
Furthermore, it will allow for the boundary terms, which are always absent in the traditional Hamiltonian picture, to be easily identified and transformed into power ports which can be used for energy exchange through the boundary of the spatial domain.

The geometric description of ideal fluid flow is then used to compile the main results of Part I of this series.
Namely, the port-Hamiltonian formulation of the kinetic energy storage subsystem and the identification of the Dirac structure that routes the energy flowing in and out of the kinetic energy  subsystem to a distributed power port (which allows the connection to subsystems that interact with the fluid within the spatial domain) and a boundary power port (which allows the connection to subsystems that allow the exchange of material).

This first paper is organized as follows:
Sec. 2 introduces the geometric description of ideal fluid flow using differential forms on the spatial domain in the quest of identifying the state space of fluid motion.
Then, Sec. 3 discusses advected quantities, their corresponding governing equations, and their own state space separated from the one of fluid motion.
Then, Sec. 4, will discuss the first two steps of the \pH modeling procedure aiming to construct the open \pH model for the kinetic-energy subsystem with variable boundary conditions and distributed stress forces.
Finally, we conclude with some remarks in Sec. 5.

In Part II, we will present the remaining steps of the \pH procedure which will utilize the distributed stress forces added to the port-Hamiltonian model presented in Sec. 4 of the present paper to add storage of internal energy for modeling compressible flow, and to add constraint forces to model incompressible flow.

\section{Kinematics of Fluid Motion}
In this section, we provide a brief introduction to the geometric description of ideal fluid flow using differential geometric tools which we will build on in this work.
For more details on this background material, the reader is referred to \citep[Ch. 8]{Abraham1988} and \citep[Ch. 11]{holm2009geometric}.

\subsection{Configuration Space and Lie Algebra}

We start by describing the physical space in which the fluid flows.
This region is mathematically represented by an $n$-dimensional compact manifold $M$, where for non relativistic physically meaningful spaces we are restricted to $n=\{2,3\}$.
The manifold could have a boundary $\partial M$ or could have no boundary ($\partial M = \emptyset$), like for instance our (almost) spherical Earth.

As the fluid flows within the domain $M$, the material particles are transported with the flow.
Consider the fluid particle passing through a point $x\in M$ in the spatial domain at time $t$, and let $v(t,x)$ denote the velocity of that particle assumed to be smooth in its arguments.
For each time $t$, $v_t := v(t,\cdot) \in \spVecM$ is a vector field on $M$.
We call $v$ the \textit{Eulerian velocity field} of the fluid.

In order to construct a configuration space that describes motion of the fluid, it is necessary to define a family of motions that move elements of the fluid along the spatial manifold $M$. 
The evolution of the fluid particle during a fixed time $t$ is described by the map $\map{g_{t}}{M}{M}$.
Thus, the particle starting at a point $x_0$ reaches the point $g_t(x_0)\in M$ after time $t$.
The map $\map{g_{t}}{M}{M}$, referred to as the \textit{flow map}, defines the current configuration of the fluid at time $t$.
The forward map $g_t:x_0\mapsto x$ takes a fluid parcel from its initial position in the reference configuration (chosen as the Lagrangian labels) to its current spatial position in the domain.
Whereas, the inverse map $\gtInv:x\mapsto x_0$ assigns the Lagrangian labels to a given spatial point.
The flow map $g_t$  is mathematically the \textit{flow} or \textit{evolution operator} of the time dependent vector field $v$ \cite[pg. 285]{Abraham1988}.

The set of all diffeomorphisms $g_t$ on $M$ is known as the \textit{diffeomorphism group} $\diffG$.
By construction the identity element $g_0$ of $\diffG$ coincides with the identity map on $M$.
The key property of $\diffG$ is that it is an (infinite-dimensional) Lie Group \cite{Marsden1970} which allows the construction of a proper analytic framework in which fluid dynamic equations can be derived in the same spirit as for finite-dimensional system.
As an example the Lie group $\diffG$ serves as the configuration space for the fluid flow in the same way as the Lie Group $SE(3)$ is the configuration space for rigid body motion.
The reader is referred to \cite{ebin1970groups} for a proper introduction to a functional-analytic arguments of considering $\diffG$ as a Lie group, which is out of the scope of this work.

The Lie algebra of $\diffG$, denoted by $\gothgV$, consists of the space of vector fields on $M$, i.e. $T_e\diffG=\spVecM$. We denote the Lie bracket of $\gothgV$ by  $[u,v]_\mathfrak{X}$, where $u,v \in \gothgV$. The algebra adjoint operator $ad_u:\gothgV \to \gothgV$ is implicitly defined by $ad_u(v):=[u,v]_\mathfrak{X}$.
The Lie bracket on $\gothgV$ is given by minus the standard Jacobi-Lie bracket on $\spVecM$, which can be expressed using the Lie Derivative as 
\begin{equation} \label{eq:bracketgothV}
[u,v]_\mathfrak{X} = ad_u(v) =  -\LieD{u}{v}.
\end{equation}

\subsection{Differential Form Representation of Lie Algebra}
\newcommand{\IsoVelToForm}{\Psi_{\normalfont\text{vol}}}

In this work we rely greatly on representing the Lie algebra $\gothgV$ using the space of differential forms on $M$.
The differential form representation is technically beneficial for a number of reasons. First, the related Grassmannian algebra operators are very efficient in carrying out calculations and proving theorems.
Second, well established results like Stokes theorem will be crucial in the development of the Dirac structures of this work. Third, the equations of motion expressed using differential forms are invariant with respect to coordinate changes.

As a Riemannian manifold, $M$ carries a metric $\bb{M}: \spVecM \times \spVecM \to \spFn{M}$ which induces a compatible volume form $\volF=*1 \in \spKForm{n}{M}$.
There are two options to identify the Lie algebra of vector fields $\gothgV = \spVecM$ with the space of $k$-differential forms $\spKForm{k}{M}$ on the Riemannian manifold $M$.
The first uses the metric $\bb{M}$ and produces the 1-form $\vf:=\bb{M}(v,\cdot) \in \spKFormM{1}$. The second option uses the interior product and the volume form to produce the ($n-1$)-form $\omega_v:=\iota_v \volF \in \spKFormM{n-1}$. For future reference, we denote the isomorphism $v \mapsto \omega_v$ by $\IsoVelToForm$.

Conversely, we will denote the vector field corresponding to an $n-1$ form $\omega \in \spKFormM{n-1}$ by $\hat{\omega} \in \spVecM$ and the vector field corresponding to a 1-form $\alpha \in \spKFormM{1}$ by $\hat{\alpha} \in \spVecM$.

The Hodge star operator makes it possible to commute between the two representations.
In particular, using the general identity
\begin{equation}\label{eq:identity_Hodge_int_product}
\iota_v \alpha = (-1)^{(k+1)(n-k)} *(\vf \wedge *\alpha),\qquad \forall v\in \spVecM, \alpha \in \spKFormM{k},
\end{equation}
and the fact that $\volF=*1 $, it holds that
\begin{equation}\label{eq:omega_v_def}
\omega_v =\iota_v \volF = \iota_v (*1) = *\vf, \qquad \text{and} \qquad \vf = \mNmOne *\omega_v.
\end{equation}

\begin{remark}
Note that the sign in (\ref{eq:identity_Hodge_int_product}) is positive if $\alpha$ is a top form or is of an odd order, which will always be the case where this identity is used in this work.
\end{remark}

\textbf{The choice we will make in this work is to represent the Lie algebra of $\diffG$, that will be denoted by $\gothg$, using $n-1$ forms.} Thus we have $\gothg = \spKFormM{n-1}$ and its corresponding Lie bracket, denoted by $\Lbrack{\cdot}{\cdot}_\gothg$, is given by the following result.

\begin{proposition}
The Lie bracket $\map{\LbrG{\cdot}{\cdot}}{\gothg\times \gothg}{\gothg}$ is given by
\begin{equation}\label{eq:Lie_bracket_g}
\LbrG{\omega}{\beta} = -\LieD{\omV}{\beta} +\divr{\omV} \beta, \qquad \omega,\beta\in\gothg,
\end{equation}
where the divergence of a vector field $v\in\spVecM$ is the function $\divr{v} \in \spFn{M}$ defined such that $\LieD{v}{\volF} = \divr{v} \volF$, and $\omV = \IsoVelToForm^{-1}(\omega) \in \gothgV$ is the vector field corresponding to $\omega \in \gothg$.
\end{proposition}
\begin{proof}
\newcommand{\betaV}{{\hat{\beta}}}
By using the map $\IsoVelToForm$ as a Lie-algebra isomorphism between $\gothgV$ and $\gothg$, the bracket on $\gothg$ can be constructed by
\begin{equation}\label{eq:LieBrackDef_g}
\Lbrack{\omega}{\beta}_\gothg := \IsoVelToForm(\Lbrack{\IsoVelToForm^{-1}(\omega)}{\IsoVelToForm^{-1}(\beta)}_\mathfrak{X}),
\end{equation}
Using (\ref{eq:bracketgothV}) and the definition of $\IsoVelToForm$, one has that
$
\LbrG{\omega}{\beta} = \IsoVelToForm([\omV,\betaV]_\mathfrak{X}) = \IsoVelToForm(-\LieD{\omV}{\betaV}) = - \iota_{(\LieD{\omV}{\betaV})} \volF.
$
Using the Lie derivative property 
$$\LieD{u}{(\iota_v \alpha)} = \iota_{(\LieD{u}{v})}\alpha + \iota_v(\LieD{u}{\alpha}), \qquad u,v\in \spVecM, \alpha\in \spKFormM{k}$$
and the definition of divergence it holds
$
\LbrG{\omega}{\beta} =-\LieD{\omV}{(\iota_\betaV \volF)} + \iota_\betaV(\LieD{\omV}{\volF})=-\LieD{\omV}{\beta} + \divr{\omV}\beta.
$
\end{proof}

\subsection{Permeable vs. Impermeable Boundaries}\label{sec:permeable_vs_imperm_boundaries}

\begin{figure}
\centering
\begin{subfigure}{0.28\textwidth}
\includegraphics[width=\textwidth]{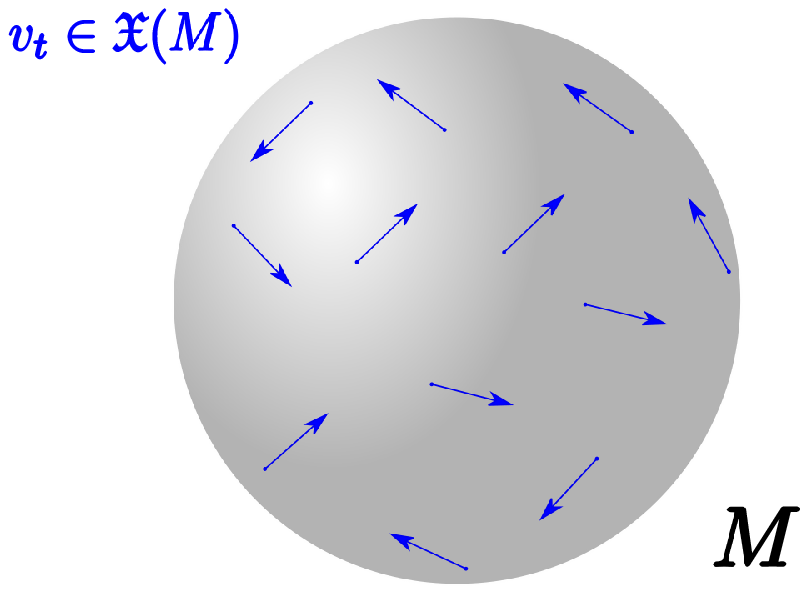}
\end{subfigure}
\qquad
\centering
\begin{subfigure}{0.28\textwidth}
\includegraphics[width=\textwidth]{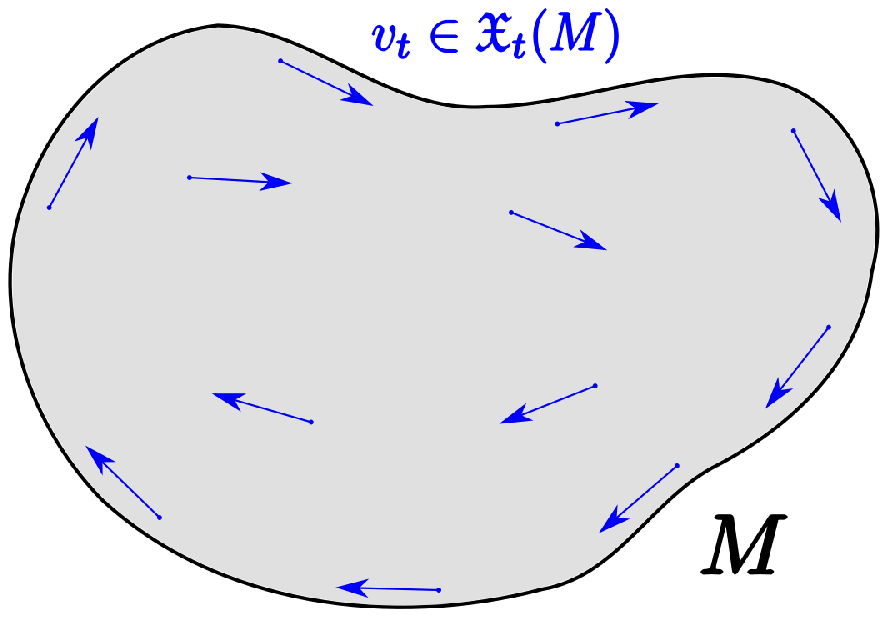}
\end{subfigure}
\qquad
\centering
\begin{subfigure}{0.28\textwidth}
\includegraphics[width=\textwidth]{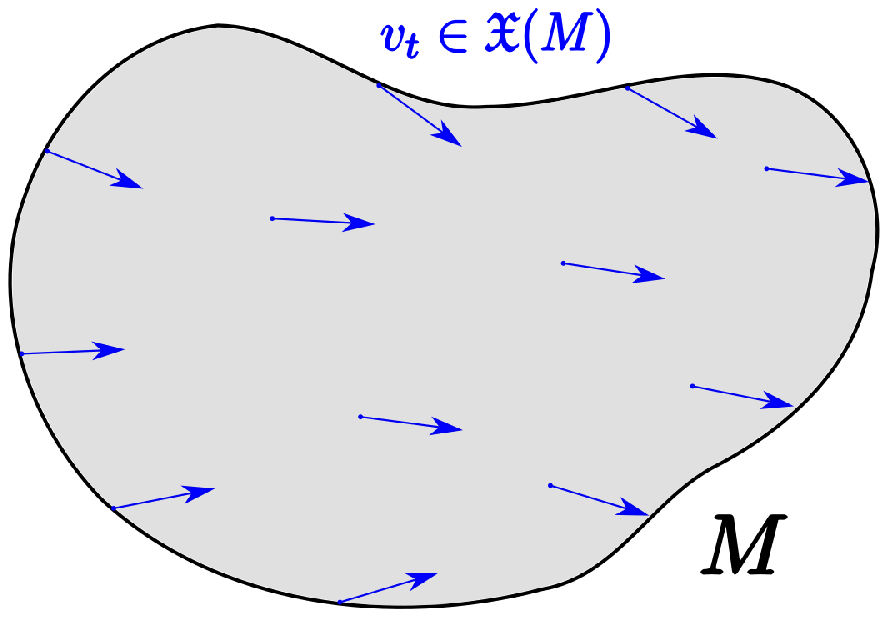}
\end{subfigure}
\caption{Boundary cases for fluid flow on a compact manifold; no boundary (left), impermeable boundary (middle), or permeable boundary (right).}
\label{fig:boundary_cases}
\end{figure}

Before moving on, we turn attention to a fundamental issue that distinguishes between the Hamiltonian and \pH treatments of fluid dynamics with respect to the boundary of the fluid container, as shown in Fig. \ref{fig:boundary_cases}.

The construction presented so far in describing the configuration of a fluid flowing using the diffeomorphism group $\diffG$ rests on the fundamental assumption that the fluid particles always remain within the fluid container.
In case $M$ has no boundary, this is obviously always true. However, in case $M$ has a boundary, the Lie algebra of $\diffG$ is constrained to the subspace $\spVecMt \subset \spVecM$ with $\spVecMt$ defined as the subset of vector fields on $M$ that are tangent to $\partial M$.
This extra condition on $v$ is synonymous with  the boundary of $M$ being impermeable.

The restriction of $\gothgV$ to vector fields that are tangent to $\partial M$ translates to the condition that its corresponding $n-1$ form is normal to $\Mbound$, i.e. satisfies $i^*(\omega_v)=0$, where $\map{i}{\Mbound}{M}$ is the canonical inclusion map.
Differential forms normal to $\Mbound$ are said to satisfy the \textit{Dirichlet boundary condition}.
Thus, in the case of impermeable boundary, the Lie algebra $\gothg$ for impermeable boundaries becomes $\spKFormMbc{n-1}{D}$, which denotes the space of normal $n-1$ forms on $M$.

The cases for which $M$ has no boundary or $\Mbound$ is impermeable correspond to an isolated fluid dynamical system that cannot exchange (mass) flow with its surrounding.
Consequently, all results of the traditional Hamiltonian formalism relying on $\diffG$ and Poisson structures are related to such isolated systems.

In order to treat fluid dynamical systems in the \pH framework, it is necessary to treat the fluid container with a boundary that is permeable by default.
This is necessary, for instance, if $M$ is to represent a control volume inside a bigger flow domain.
Fortunately, the extension of the known geometric formulation of fluids to permeable boundaries is only obstructed at the group level.
At the Lie algebra level, this extension is simply afforded by dropping the previously identified Dirichlet constraint from the vector space $\spKFormMbc{n-1}{D}$, which is the algebra for impermeable boundaries, and taking instead all of $\spKFormM{n-1}$ as the underlying vector space.

The  deceptive  simplicity  of  this  step  comes  with  a  drastic  conceptual  shift. The Lie algebra $(\gothg,\Lbrack{\cdot}{\cdot}_\gothg)$ for permeable boundaries can no longer be the associated Lie algebra of  the  diffeomorphism  group $\diffG$ or any of its subgroups.
The intuition behind it is that the flow map $g_t$ no longer becomes bijective as the fluid particles in the spatial domain $M$ comprising the reference configuration are no longer constrained to remain in $M$.
Therefore, the treatment of permeable boundaries requires nothing more and nothing less than the ability to phrase any question about the fluid in a form that can be treated entirely at Lie algebra level.

The discussion of impermeable and permeable boundaries at this early stage is important, because much of what follows depends on whether one deals with impermeable or permeable boundaries.
For instance, the determination of the duals of the vector space underlying the Lie algebra $\gothg$, crucially depends on the nature of $\Mbound$. Another consequence is that boundary terms, which all vanish for impermeable boundaries, must be carried through in a treatment that aspires to be valid for the three boundary cases. 

Indeed, it is precisely the surface terms that make the key difference between the Lie-Poisson structure, underlying the traditional Hamiltonian theory, and the Stokes-Dirac structure, underlying the \pH theory, and will be of extreme importance in this work as it will rigorously reveal the corresponding boundary ports.

\subsection{Dual Space of Lie Algebra}
The dual space $\gothgstar$ to the Lie algebra $\gothg$ is naturally identified with the space of 1-forms $\spKFormM{1}$ by means of the duality pairing
\begin{equation}
\pair{\alpha}{\omega}_\gothg := \int_M \alpha \wedge \omega, \qquad \alpha \in \gothgstar, \omega \in \gothg,
\end{equation}
which is the natural pairing based on integration on $M$ of differential forms of complementary order.

An essential ingredient for deriving the equations of fluid motion in the Hamiltonian formalism is the map $\map{ad^*_\omega}{\gothgstar}{\gothgstar}$, which is the formal dual to the adjoint operator of $\gothg$ defined by $ad_{\omega}(\cdot) := \LbrG{\omega}{\cdot}$, for a given $\omega \in \gothg$.
The explicit expression for $ad^*_\omega$ is given by the following result. 

\newcommand{\surfAdStar}[1]{\eta_{ad_{#1}}}

\begin{proposition}\label{prop:adStar}
For any $\omega \in \gothg$, the dual map $ad_\omega^*$ of the adjoint operator $ad_\omega$ of $\gothg = \spKForm{n-1}{M}$ is given by
\begin{equation}\label{eq:adStar}
ad^*_{\omega}(\alpha) = \LieD{\hat{\omega}}{\alpha} + \divr{\hat{\omega}} \alpha ,\qquad \omega \in \gothg, \alpha \in \gothgstar.
\end{equation}
For a general manifold $M$ with boundary $\partial M$, the map $ad_\omega^*$ in (\ref{eq:adStar}) satisfies for any $\omega, \beta \in \gothg$ and $ \alpha \in \gothgstar,$ 
\begin{equation}\label{eq:adStar_pairing}
\pair{ad^*_{\omega}(\alpha)}{\beta}_\gothg =  \pair{\alpha}{ad_{\omega}(\beta)}_\gothg + \int_{\Mbound} \bound{\surfAdStar{\omega}(\alpha,\beta)},\qquad \omega, \beta \in \gothg, \alpha \in \gothgstar,
\end{equation}
where the $(n-1)$-form $\surfAdStar{\omega}(\alpha,\beta) \in \spKFormM{n-1}$ is given by
\begin{equation}\label{eq:surfAd}
\surfAdStar{\omega}(\alpha,\beta) = \iota_{\hat{\omega}}(\alpha\wedge\beta).
\end{equation}
The term $\bound{\eta} \in \spKForm{k}{\Mbound}$ denotes the trace of the form $\eta\in\spKFormM{k}$ which is defined as the pullback of the inclusion map $\map{i}{\Mbound}{M}$, i.e. $\bound{\eta}:= i^*\eta $.

\end{proposition}
\begin{proof}
\newcommand{\omB}{\beta}
Consider any $\omB \in \gothg$. By using (\ref{eq:Lie_bracket_g}) we have that
\begin{align}
\pair{\alpha}{ad_{\omega}(\omB)}_\gothg = \int_M \alpha \wedge \LbrG{\omega}{\omB}
										= \int_M \alpha \wedge \divr{\omV} \omB- \int_M \alpha \wedge \LieD{\omV}{\omB}.\label{eq:ad_proof_1}
\end{align}
Using the Leibniz rule for the Lie derivative
$
\LieD{\omV}{(\alpha \wedge \beta)} = \LieD{\omV}{\alpha}\wedge \beta + \alpha \wedge \LieD{\omV}{\beta}$
and Cartan's magic formula $\LieD{\omV}{}=\extd \iota_{\omV} +\iota_{\omV} \extd$ the second integrand in (\ref{eq:ad_proof_1}) becomes
\begin{align}
- \alpha \wedge \LieD{\omV}{\omB}	= \LieD{\omV}({\alpha} \wedge \omB) - \extd \iota_{\omV}(\alpha \wedge \omB).\label{eq:ad_proof_3}
\end{align}
Substituting (\ref{eq:ad_proof_3}) in (\ref{eq:ad_proof_1}), using the definition of divergence and Stokes theorem we have that
\begin{align*}
\pair{\alpha}{ad_{\omega}(\omB)}_\gothg =& \pair{\divr{\omV} \alpha +\LieD{\omV}{\alpha}}{\omB}_\gothg - \int_{\partial M} i^* (\iota_{\omV}(\alpha \wedge \omB))\\
										=& \pair{ad^*_{\omega}(\alpha)}{\omB}_\gothg - \int_{\partial M} \bound{(\iota_{\omV}(\alpha \wedge \omB))}.\label{eq:ad_proof_4}
\end{align*}
\end{proof}

\begin{remark}
Interestingly (but not surprisingly), this condition on representatives of the Lie algebra in case of impermeable boundary, is exactly the one that nullifies the boundary term in (\ref{eq:adStar_pairing}), which is easily verified because the pullback in the trace distributes over the wedge. This consideration sheds light on the fact that in classical Hamiltonian theory the state space is constrained such that no energy exchange can happen on the boundary of the spatial manifold, and the surface term in the duality pairing (\ref{eq:adStar_pairing}) is consequently neglected. 

As will be shown in what follows, the surface terms are the key difference between the Lie-Poisson structure, underlying the traditional Hamiltonian theory, and the Stokes-Dirac structure, underlying the port-Hamiltonian theory, and will be of extreme importance in this work. Our methodology does not rely on constraining the state space of the fluid, but to consider duality pairings like (\ref{eq:adStar_pairing}) in their full generality, which will reveal rigorously boundary ports that can be used for modeling or control purposes. 

\end{remark}

\section{Kinematics of Advected Quantities}
In ideal continuum flow, the material parcels of the fluid are carried by the flow.
These parcels are transported (advected) by the ideal flow along with extensive thermodynamic properties such as the parcels' mass and heat.
{The properties transported by the flow are referred to as \textbf{advected quantities}.

The presence of advected quantities adds more structure to the spaces underlying the motion of fluids $\gothg$ and $\gothgstar$.
Understanding this additional structure is crucial for the development of the decomposed model of fluid dynamics we aim for in this work as well as the derivation of the Stokes-Dirac structure underlying the kinetic energy subsystem.

Next, we describe mathematically advected quantities and their governing evolution equations.
Then, we introduce the interconnection maps that will allow relating the spaces of fluid motion and the spaces of advected quantities to each other.
Finally, we discuss the added \textit{semi-direct algebra} structure that the presence of these advected quantities introduces into the geometric picture.
%

\subsection{Mathematical Description}\label{sec:math_desc_advection}

In general, an advected quantity $a_t$ is represented mathematically as a time-dependent tensor field. We denote the vector space of advected quantities by $V^* \subset \cl{T}(M)$ which is usually a subspace of the space of tensor fields $\cl{T}(M)$ on $M$.
Specific examples include scalar fields (e.g. buoyancy, entropy), vector fields (e.g. magnetic field), 2-forms (e.g. vorticity), and top forms (e.g. mass form).
All advected quantities that will be considered in this work will be represented by differential forms.

The relation between the differential form $a_t$ at $t>0$ and its initial value $a_0$ is given by the pullback of the flow map $g_t \in G$.
Thus, we have that $a_0 = g_t^* a_t$ or $a_t = (g_t^{-1})^* a_0$.
An important observation is that the Eulerian advected quantity $a_t$ is thus completely determined by the flow map $g_t$ and its initial value $a_0$.

The differential expression for $a_t$ to be advected with the flow is that
\begin{equation}\label{eq:advQuantity}
(\frac{\partial}{\partial t} + \LieD{v}{}) a_t = 0,
\end{equation}
which is an immediate consequence of the Lie derivative formula for time-dependent tensor fields \citep[Pg. 372]{Abraham1988}.
An advected quantity is also referred to as being invariant under the flow, frozen into the fluid, or  Lie dragged with the flow.

The prototypical case of an advected quantity that occurs in all continuum flows is the mass top-form given by $\mForm := \rho_t \volF \in \spKForm{n}{M}$, where $\rho_t \in \spFn{M}$ denotes the mass density function.
We have that $\mForm \in V^* = \spKForm{n}{M}$.

The mass conservation (or continuity) equation  can be written in terms of the mass form as $g_t^* \mForm = \mu_0 $, or equivalently 
\begin{equation}\label{eq:cnty_mass_form}
\pt \mForm + \LieD{v}{\mForm} = 0.
\end{equation}

To express the mass continuity in terms of the density function $\rho_t$, we substitute $\mForm = \rho_t \volF$ in (\ref{eq:cnty_mass_form}) which can be manipulated to result in
\begin{equation}\label{eq:continuity_density}
\pt \rho_t + \LieD{v}{\rho_t} + \rho_t \divr{v} = 0.
\end{equation}

By comparing the mass continuity equation in its two forms (\ref{eq:cnty_mass_form}) and (\ref{eq:continuity_density}), it is observed that in general, the mass form $\mForm$ is an advected quantity while the mass density function $\rho_t$ is not.
This observation is the main reason why $\mForm$ will be chosen as a state variable later in this work for representing the fluid's kinetic energy, which depends on the density and velocity of the flow.

Other examples of advected quantities, which will be used in Part II, include the entropy function $s_t\in \spFn{M}$, in case of adiabatic compressible flow, and the volume form $\volF \in \spKFormM{n}$ in case of incompressible flow.

\subsection{Relation Between Fluid Motion and Advected Quantities}\label{sec:interconn_maps}

\newcommand{\mapR}{\tilde{\varphi}_a}
\newcommand{\mapRdual}{\tilde{\varphi}_a^*}

\begin{figure}
\centering
\includegraphics[width=0.75\textwidth]{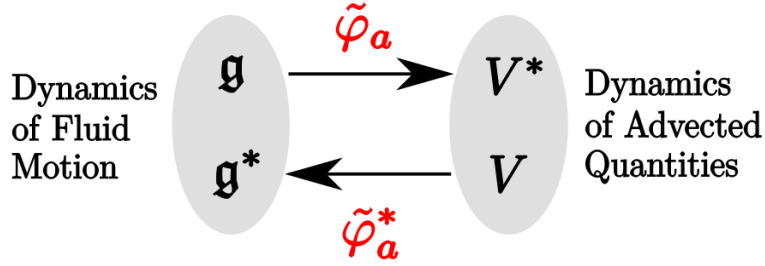}
\caption{Interconnection maps relating the space of fluid motion to the space of advected quantities.}
\label{fig:semidirect_Maps}
\end{figure}

In order to build a decomposed model of fluid flow, it is important to describe the dynamics of advected quantities separately from the dynamics of fluid motion.
With reference to Fig. \ref{fig:semidirect_Maps}, the dynamics of the fluid motion are defined on the space $(\gothg\times\gothgstar)$, whereas the dynamics of advected quantities are defined on $(V^*\times V)$.
The dual space of $V^* = \spKFormM{k}$ is $V = \spKFormM{n-k}$ with respect to the duality pairing $\map{\pair{\cdot}{\cdot}_{V^*}}{V\times V^*}{\bb{R}}$ given by the integral of the wedge product, i.e.
\begin{equation}
\pair{\bar{a}}{a}_{V^*} := \int_M \bar{a} \wedge a, \qquad a\in V^*, \bar{a}\in V.
\end{equation}

As discussed in Sec. \ref{sec:math_desc_advection}, the evolution of an advected quantity depends on the fluid motion described by $\omega\in\gothg$.
On the other hand, the advected quantity carried by the flow influences the flow motion in a power-consistent manner through the \textit{bidirectional exchange} of kinetic energy of the fluid and the additional energy characterized by the advected quantity (e.g. potential or magnetic energy).

The effect of the fluid motion ($\omega\in\gothg$ on the advected quantity $a\in V^*$ is encoded in the the primary map $\map{\mapR}{\gothg}{V^*}$, which is defined by
\begin{equation}\label{eq:def_Rmaps}
\fullmap{\mapR}{\gothg}{V^*}{\omega}{\mapR(\omega) := \LieD{\omV}{a}.}
\end{equation}
In terms of this primary map, the governing equation of an advected quantity (\ref{eq:advQuantity}) is rewritten as $\dot{a}_t = - \mapR(\omega)$, with $\omega = \IsoVelToForm(v)$.

The reverse effect, which the advected quantities has on the fluid motion, is characterized by the (formal) dual map of $\mapR$ given by
\begin{equation}\label{eq:def_R_dual_maps}
\fullmap{\mapRdual}{V}{\gothgstar}{\abar}{\mapRdual(\abar).}
\end{equation}

We refer to the two maps $\mapR,\mapRdual$ as the \textbf{interconnection maps} (\textit{cf.} Fig. \ref{fig:semidirect_Maps}) as they will serve the fundamental role of interconnecting the space of advected quantities with the space of fluid motion in the port-Hamiltonian model, as will be shown later in this part and in Part II.

%




\newcommand{\surfPhiDTA}[1]{\eta_{\tilde{\varphi}_{#1}}}

The explicit expression for $\mapRdual$ is given by the following result which depends on a specific choice of $V^*$.
For our purpose, we only consider the cases of top-forms and smooth functions relevant to the advected quantities of mass form ($V^* = \spKFormM{n}\ni \mu $) and entropy ($V^* = \spKFormM{0}\ni s $).
\begin{proposition}\label{prop:diamond}
Let $V$ be a representation space of $\diffG$ and $V^*$ its dual space, where both $V$ and $V^*$ are spaces of differential forms.
For a given $a \in V^*$, consider the map $\mapR$ defined by (\ref{eq:def_Rmaps}) and consider its dual map  $\map{\mapRdual}{V}{\gothgstar}$.
For a general manifold $M$ with boundary $\partial M$, the dual map  $\mapRdual$ satisfies for any $\omega\in\gothg$
\begin{equation}
\pair{\mapRdual(\abar)}{\omega}_\gothg = \pair{\abar}{\mapR(\omega)}_{V^*}+ \int_{\Mbound} \bound{\surfPhiDTA{a}(\omega,\abar)},\label{eq:pair_Rmap}
\end{equation}
where $\surfPhiDTA{a}(\omega,\abar) \in \spKFormM{n-1}$ is the corresponding $n-1$ form representing the surface term.

\begin{enumerate}
\item In case $V = \spKFormM{0}$ and $V^* = \spKFormM{n}$, 
\begin{equation*}
\mapRdual(\abar) = -(*a) \extd \abar, \qquad   \surfPhiDTA{a}(\omega,\abar) = -(*a)\omega\wedge \abar,
\end{equation*}
\item In case $V = \spKFormM{n}$ and $V^* = \spKFormM{0}$, 
\begin{equation*}
\mapRdual(\abar) =  (*\abar) \extd a, \qquad   \surfPhiDTA{a}(\omega,\abar) = 0,
\end{equation*}

\end{enumerate}
where $*$ denotes the Hodge star operator and $\extd$ denotes the exterior derivative operator.
\end{proposition}
\begin{proof}
\begin{enumerate}
\item Let $\abar \in V = \spKFormM{0}$ and $a\in V^*=\spKFormM{n}$, where $a$ can be written as $a= (*a)\volF$.
Let $\omega\in \gothg$ and let $\omV \in \spVecM$ be its corresponding vector field.

Using the fact that $\abar\in \spFn{M}$ and the Leibniz rule for the exterior derivative, we can write 
\begin{align*}
\pair{\abar}{\mapR(\omega)}_{V^*} &= \int_M \abar\wedge \LieD{\omV}{a} = \int_M \abar \wedge \extd \iota_{\omV}{a}= \int_M \abar \wedge \extd ((*a) \iota_{\omV}{\volF}),\\
									&= \int_M \abar \wedge \extd (*a \omega) = \int_{M} \extd(\abar\wedge (*a \omega) ) -   \extd \abar\wedge( *a) \omega ,\\
									&= \int_{\Mbound} (*a) \omega\wedge \abar - \int_{M}  (*a) \extd \abar \wedge \omega ,
									= \int_{\Mbound} - \surfPhiDTA{a}(\omega,\abar) +  \pair{ - (*a) \extd \abar}{\omega}_\gothg ,
\end{align*}
which concludes the proof of (i).

\newcommand{\bbar}{{\bar{b}}}
\item Let $\bbar \in V = \spKFormM{n}$ and $b\in V^*=\spKFormM{0}$, where $\bbar$ can be written as $\bbar= (*\bbar)\volF$.

Using the definition of $\tilde{\varphi}_b(\omega)$ and the fact that $b\in \spFn{M}$, we write
\begin{equation*}
\pair{\bbar}{\tilde{\varphi}_b(\omega)}_{V^*} = \int_M \bbar \wedge \LieD{\omV}{b} =  \int_M \bbar \wedge \iota_\omV \extd b= \int_M \iota_\omV \extd b \wedge \bbar.
\end{equation*}
Using the Leibniz rule for the interior product it holds
$\iota_\omV(\extd b \wedge \bbar) = \iota_\omV(\extd b) \wedge \bbar - \extd b \wedge \iota_\omV(\bbar)=0,$
since, as an $n+1$ form, $\extd b \wedge \bbar=0$. It follows
\begin{equation*}
    \int_M \iota_\omV \extd b \wedge \bbar =  \int_M \extd b \wedge \iota_\omV \bbar	=  \int_M \extd b \wedge (*\bbar) \iota_\omV \volF 
    = \int_M (*\bbar) \extd b \wedge  \omega=\pair{\tilde{\varphi}_b^*(\bbar)}{\omega}_\gothg
\end{equation*}
Moreover, from (\ref{eq:pair_Rmap}) we have that $\surfPhiDTA{b}(\omega,\bbar) = 0$.
\end{enumerate}
\end{proof}

\subsection{Semidirect Product Structure}
\newcommand{\elS}[1]{(\omega_{#1},\abar_{#1})}
\newcommand{\elSdual}[1]{(\alpha{#1},a_{#1})}
\newcommand{\surfLie}[3]{\iota_{#1}(#2\wedge #3)}
\newcommand{\surfAdSstar}[1]{\eta_{\boldsymbol{ad}_{#1}}}

The presence of advected quantities adds more structure to the Lie algebra structure of $\gothg$ that will be central for deriving the Hamiltonian dynamics in the coming section.
On the group level (valid only for impermeable or no boundary), the pullback operation of a flow map $g_t \in \diffG$ on an advected quantity $a_t\in V^*$ defines a right linear action (i.e. representation) of the group $\diffG$ on the vector space $V^*$, which induces another representation on its dual space $V$, given also by the pullback operation \citep{holm2009geometric}.
Thus, both $V$ and $V^*$ are representation spaces of $\diffG$.

On the algebra level, the representation of $\diffG$ on $V$ and $V^*$ induces two other representations of $\gothg$ on  $V$ and $V^*$ that are both related to the Lie derivative operator \citep{holm2009geometric}.  

Given the group $\diffG$, vector space $V = \spKFormM{k}$, and the right representation given by the pullback operation, we define the semi-direct product Lie group $S$ as the group with underlying manifold $\diffG\times V$ and group operation $\map{\bullet}{S\times S}{S}$ defined by
\begin{equation}
(g_1,\abar_1)\bullet (g_2,\abar_2) = (g_1\circ g_2, {g_2}^*(\abar_1) + \abar_2).
\end{equation}
Usually $S$ is denoted by $\diffG\ltimes V$.
The Lie algebra of $S$ is then given by the semi-direct product algebra $\goths = \gothg \ltimes V$, with its corresponding bracket $\map{\LbrS{\cdot}{\cdot}}{\goths \times \goths}{\goths}$, is defined using the Lie bracket on $\gothg$ in (\ref{eq:Lie_bracket_g}) and the induced action of $\gothg$ on $V$ as
\begin{equation}\label{eq:Liebracket_s}
\LbrS{\elS{1}}{\elS{2}} := (\LbrG{\omega_1}{\omega_2},\LieD{\hat{\omega}_2}{\abar_1}-\LieD{\hat{\omega}_1}{\abar_2}).
\end{equation}

The dual space $\gothsStar = \gothgstar\times V^*$ to the Lie algebra $\goths$ will serve as the Poisson manifold on which the dynamical equations of motion will be formulated later.
The duality pairing between an element $(\alpha,a)\in \gothsStar$ and an element $\elS{} \in \goths$ is given by
\begin{equation}\label{eq:duality_s}
\pair{\elSdual{}}{\elS{}}_\goths := \pair{\alpha}{\omega}_\gothg + \pair{a}{\abar}_V = \int_M \alpha \wedge \omega + a \wedge \abar,
\end{equation}
where the duality pairing $\map{\pair{\cdot}{\cdot}_V}{V^*\times V}{\bb{R}}$ is given by
\begin{equation}
\pair{a}{\bar{a}}_V := \int_M a \wedge \bar{a}, \qquad a\in V^*, \bar{a}\in V.
\end{equation}

As usual we define the adjoint operator $\map{\adS{\elS{}}}{\goths}{\goths}$ for a given $\elS{} \in \goths$ by
$
\adS{\elS{}} := \LbrS{\elS{}}{\cdot}.
$
Note that the \textit{bold} notation for $\adS{\elS{}}$ is to distinguish it from the $ad_\omega$ operator of the Lie algebra $\gothg$.
Then the dual map $\map{\adSdual{\elS{}}}{\gothsStar}{\gothsStar}$ is given by the following result, which is an extension of Prop. \ref{prop:adStar} for the case of the semi-direct product $\goths$.
\begin{theorem}\label{theorem:adStarS_duality}
For a given $\elS{} \in \goths$, the formal dual $\adSdual{\elS{}}$ of the adjoint operator $\adS{\elS{}}$ and (\ref{eq:Liebracket_s}) is given explicitly by
\begin{equation}\label{eq:adstar_s}
\adSdual{\elS{}}\elSdual{} = (ad^*_{\omega}(\alpha) + \abar \diamond a, \LieD{\hat{\omega}}{a}), \qquad \elSdual{}\in \gothsStar,
\end{equation}
where $\map{ad^*_{\omega}}{\gothgstar}{\gothgstar}$ is given by Prop. \ref{prop:adStar}, and the diamond map $\map{\diamond}{V\times V^*}{\gothgstar}$ is given by
$$\abar \diamond a := (-1)^{c+1} \mapRdual(a), \qquad \abar \in V = \spKFormM{k}, a\in V^*= \spKFormM{n-k},$$
where $c= k(n-k) \in \bb{R}.$

For a general n-dimensional manifold $M$ with boundary $\partial M$, the map $\adSdual{\elS{1}}$ satisfies the following equality for any $\elS{1},\elS{2}\in \goths$ and $\elSdual{}\in \gothsStar$
\begin{equation}\label{eq:pairing_adSdual}
\pair{\adSdual{\elS{1}} \elSdual{}}{\elS{2}}_\goths  = \pair{\elSdual{}}{\adS{\elS{1}}\elS{2}}_\goths + \int_{\Mbound}\bound{ \surfAdSstar{\elS{1}}(\alpha,a,\omega_2,\abar_2)},
\end{equation}
where the surface term $\surfAdSstar{\elS{1}}(\alpha,a,\omega_2,\abar_2) \in \spKFormM{n-1}$ is expressed as
\begin{equation}\label{eq:surfAdS}
\surfAdSstar{\elS{1}}(\alpha,a,\omega_2,\abar_2) = \surfAdStar{\omega_1}(\alpha,\omega_2) +(-1)^{c+1}  \surfPhiDTA{a}(\omega_2,\abar_1) - \surfLie{\omega_2}{a}{\abar_1}  + \surfLie{\omega_1}{a}{\abar_2} ,
\end{equation}
where the expressions of the surface terms are given in (\ref{eq:surfAd}) and in Prop. \ref{prop:diamond}. 
\end{theorem}
\begin{proof}
See Appendix Sec. \ref{proof:adStarS_duality}.

\end{proof}

\begin{remark}
The diamond operator $\map{\diamond}{V\times V^*}{\gothgstar}$ in (\ref{eq:adstar_s}) is usually introduced in the Hamiltonian mechanics literature, e.g. \citep{marsden1984semidirect,holm1998euler,Modin2011}, as a short-hand notation for the dual map (\ref{eq:def_R_dual_maps}) and it describes the effect that the advected quantities impose on the motion of the fluid.

In this notation of the diamond operator, the element of $V^*$ plays the role of modulation in the linear maps  (\ref{eq:def_Rmaps}) and (\ref{eq:def_R_dual_maps}).
Furthermore, the image of $\omega$ under $\mapR$ belongs to the tangent space at $a\in V^*$, i.e.
$\mapR(\omega) = -\dot{a} \in T_a V^* \cong V^*.$
On the other hand, the range of the map $\mapRdual$ is the cotangent space at $a\in V^*$, i.e. $T_a^*V^* \cong V$.

\end{remark}


\section{Port-Hamiltonian Modeling of the Kinetic Energy Subsystem}
Now we start the \pH modeling procedure of a fluid dynamical system flowing on an arbitrary spatial manifold.
Unlike the standard \textit{top-down} approach in classical Hamiltonian theory for modeling a fluid system starting from the configuration space $\diffG$, the philosophy of the port-Hamiltonian modeling process is based on a \textit{bottom-up} approach.

The \pH procedure for modeling a general dynamical system can be summarized as follows:

\begin{enumerate}
\item \textbf{Conceptual tearing:} The starting point for the \pH approach for modeling is a \textit{conceptual tearing} process of the overall physical system, viewing it as a set of interconnected energetic subsystems. The $i$-th subsystem is characterized by the energy it possess, denoted by  $\map{H_i}{\cl{X}_i}{\bb{R}}$. This energy traverses from one subsystem to another through an imaginary boundary dividing both of them.
 
\item \textbf{Hamiltonian modeling of isolated energetic subsystems:} After identifying the individual energetic subsystems, one can use the standard Hamiltonian theory at this point to develop a Hamiltonian model for each energetic subsystem \textit{isolated} from the rest of the system. However, the emphasis here is on physical energy variables $x_i\in \cl{X}_i$ which do not necessarily correspond to the canonical coordinates of some cotangent bundle.
In fact, energy variables that are most physically compelling are in general non-canonical \citep{Morrison1998}.
Consequently, the outcome of this stage is a closed Hamiltonian model for each energetic subsystem, defined by the energy function $H_i(x_i)$ and its corresponding non-canonical Poisson structure.

\item \textbf{Add interaction ports:} The closed Hamiltonian model of each subsystem is now extended to an open \pH model with interaction ports. One then replaces the Poisson structure representing the conservation of $H_i$ by a Dirac structure that allows for non-zero energy exchange via the interaction ports.

\item \textbf{Interconnect all energetic subsystems:} The overall physical system is then constructed by interconnecting the different energetic subsystems. This is achieved by specifying the interconnection structure of the open \pH subsystems in the form of constraints on the variables of the interaction ports of each subsystem. For the interconnection to be power-consistent, only subsystems with \textit{compatible} interaction ports are interconnected, where the geometric structure of Lie algebras plays a significant role. Additional Dirac structures could be possibly used to resolve incompatibility issues.

\item \textbf{Compact the \pH model:} After developing a decomposed model of the physical system in terms of a network of interconnected subsystems, one could (optionally) compact the network by combining the energy storing elements together and all the Dirac structures into one, denoted by $D(x)$.
The compact \pH model would then have a total state space $\cl{X}$ and total Hamiltonian $\map{H}{\cl{X}}{\bb{R}}$ given, respectively, by
$$\cl{X} := \prod_i \cl{X}_i, \qquad H(x):= \sum_{i} H_i(x_i), $$
where $x \in \cl{X}$ is the state of the overall system corresponding to the combined energy variables, which may have redundancies.
The compact model could be either a \textit{closed} system, or an \textit{open} one with interaction ports that allows the overall \pH model to be interconnected to its external environment.
\end{enumerate}

For the fluid dynamical system at hand, the conceptual tearing process yields two energetic subsystems; One for storage of kinetic energy and another for storage of internal energy (for the general case of compressible flow).
In what remains of Part I, we will only focus on the application of the second and third steps for modeling the kinetic energy subsystem.
Indeed, we will derive it as a closed Hamiltonian model and then extend it to an open \pH model.

The non-canonical coordinates we choose to develop the Hamiltonian model belong to the dual space $\gothsStar = \gothgstar\times V^*$ of the semi-direct Lie-algebra $\goths=\gothg\ltimes V$.
In the classical fluid dynamics literature, these coordinates are referred to as the Eulerian representation of the fluid motion, in contrast to the Lagrangian representation corresponding to the canonical coordinates on the cotangent bundle of $\diffG$.

Using the \textit{semi-direct product reduction theorem} \citep{marsden1984semidirect,marsden1984reduction}, we will show in the coming section that the Hamiltonian dynamical equations in terms of the Eulerian variables is derivable from the canonical symplectic Hamiltonian dynamics in terms of the Lagrangian variables. 
This procedure will yield the Hamiltonian dynamics on the dual space $\gothsStar = \gothgstar\times V^*$ of the semi-direct Lie-algebra $\goths=\gothg\ltimes V$.
It is important to note that this is the only stage where the diffeomorphism group is used, which is possible since we will deal with an isolated closed dynamical system.
Moreover, the material presented in this work will be a mere reproduction of the results of \citep{marsden1984semidirect,marsden1984reduction} with a slight difference that the algebra considered here is $n-1$ differential forms.

Afterwards, we will show the derived Hamiltonian model and its corresponding Lie-Poisson structure are extended to an open \pH model based on a Stokes-Dirac structure. In this open model, the group structure will be no longer valid, but all the constructions on the algebra will be extended directly to incorporate the permeable boundaries.

\subsection{Closed Model of Kinetic Energy}


The \textit{kinetic co-energy}\footnote{For the motivation behind using the co-energy and energy terminologies, see \cite[Sec. B.2]{duindam2009modeling}.} of the fluid flow can be used to construct a (reduced) Hamiltonian functional on $\gothsStar = \gothgstar\times V^*$.
By analogy with the kinetic co-energy of a system of particles, we define the co-energy of the fluid flowing on a manifold $M$ to be
\begin{equation}
E^*_k(v,\rho) := \int_M \half \rho \bb{M}(v,v)\volF,
\end{equation}
where $\bb{M}$ is the Riemannian metric on $M$.

Now we wish to represent the kinetic co-energy of the fluid in terms of the Lie algebra $\gothg$ and the space of advected quantities $V^*$, which is expressed as
\begin{equation}\label{eq:reduced_Lagrangian}
L_k(\omega_v,\mu) = \int_M \half (*\mu) \omega_v \wedge * \omega_v,
\end{equation}
which follows from $\rho = *\mu$ and the equalities
\begin{equation}\label{eq:equalities_kinetic_energy}
\bb{M}(v,v)\volF = \iota_v \vf \volF =  \vf \wedge *\vf = \mNmOne *\omega_v \wedge \omega_v = \omega_v \wedge * \omega_v,
\end{equation}
where the second equality follows from identity (\ref{eq:identity_Hodge_int_product}).
Thus, $\map{L_k}{\gothg\times V^*}{\bb{R}}$ is a (reduced Lagrangian) functional on the Lie algebra of $\diffG$ that depends parametrically on $\mu \in V^*$.
Using a partial Legendre transformation, we can thus of course represent the \textit{kinetic energy} as a functional on $\gothgstar\times V^*$ as follows.

\begin{theorem}\label{theorem:Hamiltonian}
The Legendre transformation of $\map{L_k}{\gothg\times V^*}{\bb{R}}$, as given by (\ref{eq:reduced_Lagrangian}), is the functional $\map{H_k}{\gothgstar\times V^*}{\bb{R}}$ given by
\begin{equation}\label{eq:reduced_Hamiltonian}
H_k(\alpha,\mu) = \int_M \frac{1}{2(*\mu)} \alpha \wedge * \alpha,
\end{equation}
where $\alpha:= \mNmOne (*\mu)*\omega_v \in \gothgstar$ denotes the momentum of the fluid.

The variational derivatives $\delta_\alpha H_k \in T_\alpha^*\gothgstar \cong \gothg = \spKFormM{n-1}$ and $\delta_\mu H_k \in T_\mu^*V^* \cong V =  \spKFormM{0}$ with respect to the states $\alpha \in \gothgstar $ and $\mu \in V^*$, respectively, are given by
\begin{equation}\label{eq:var_der_H_k}
\delta_\alpha H_k = \frac{*\alpha}{*\mu}, \qquad\qquad \delta_\mu H_k = - \frac{1}{2 (*\mu)^2} \iota_{\hat{\alpha}} \alpha,
\end{equation}
where $\hat{\alpha} \in \spVecM$ denotes the vector field corresponding to the 1-form $\alpha$.
\end{theorem}
\begin{proof}
See Appendix Sec. \ref{proof:Hamiltonian}.
\end{proof}

One of the benefits of the previous systematic construction of the (reduced) Hamiltonian $H_k$ is that the correct \textit{conjugate momentum} variable with respect to $\omega_v \in \gothg$ is chosen, and all variational derivatives are correctly derived.
As for the governing equations of motion, the semi-direct product reduction theorem asserts that the Hamiltonian dynamics on $\gothsStar = \gothgstar\times V^*$ is given as follows.

\begin{theorem}\label{theorem:Ham_dynamics}
The Hamiltonian dynamical equations of fluid flow in terms of Eulerian variables on the dual space $\gothsStar = \gothgstar\times V^*$ are given by
\begin{equation}\label{eq:Ham_kin_subs}
\TwoVec{\dot{\alpha}}{\dot{\mu}} = - \adSdual{(\delta_\alpha H_k,\delta_\mu H_k )}\state = \TwoVec{-ad^*_{\delta_\alpha H_k}(\alpha) - \delta_\mu H_k \diamond \mu}{-\LieD{(\delta_\alpha H_k)^\wedge}{(\mu)}},
\end{equation}
where the Hamiltonian $H_k\state$ and its corresponding variational derivatives are given by Theorem \ref{theorem:Hamiltonian}. 

The corresponding Lie-Poisson bracket underlying these dynamical equations is given by
\begin{equation}\label{eq:Lie-Poisson_bracket_kin}
\Pbrack{F}{G} (\alpha,\mu) = \pair{\alpha}{\LbrG{\delta_\alpha F}{\delta_\alpha G}}_\gothg + \pair{\mu}{\LieD{(\delta_\alpha G)^\wedge}{\delta_\mu F} - \LieD{(\delta_\alpha F)^\wedge}{\delta_\mu G}}_V,
\end{equation}
where $F,G \in \spFn{\gothsStar}$ are smooth functions on $\gothsStar$.
\end{theorem}
\begin{proof}
To derive the Hamiltonian dynamics using the semidirect product reduction theorem, we need the Hamiltonian $\B{H}_{a_0}$ described on the phase space $T^*\diffG$.
The goal is to extend the Hamiltonian functional $\map{H_k}{\gothgstar\times V^*}{\bb{R}}$ in (\ref{eq:reduced_Hamiltonian}) to one on $T^*\diffG$.
Consider the isomorphism $\map{\Gamma}{T^*\diffG\times V^*}{\gothgstar\times V^*}$ defined by
\begin{equation}
\Gamma(g,\pi,\mu_0) = ((R_g)^* (\pi), (g^{-1})^* \mu_0),
\end{equation}
which physically represents the map from the space ($T^*\diffG\times V^*$) of Lagrangian coordinates to the space ($\gothgstar\times V^*$) of Eulerian coordinates.
The action of $\Gamma$ on $(g,\pi)\in T^*\diffG$ is given by the pullback of the right translation map $\map{R_g}{\diffG}{\diffG}$, whereas the action of $\Gamma$ on $\mu_0 \in V^*$ is given by the pullback operator.
Then, by construction, the functional $\B{H}(g,\pi,\mu_0):= H_k \circ \Gamma (g,\pi,\mu_0) = H_k(\alpha,\mu)$ is a Hamiltonian on $T^*\diffG\times V^*$.
Moreover, define the Hamiltonian $\map{\B{H}_{a_0}}{T^*\diffG}{\bb{R}}$ by $\B{H}_{\mu_0}(g,\pi) = \B{H}(g,\pi,\mu_0)$.
By construction, $\B{H}_{a_0}$ is right-invariant under the action of the stabilizer group $G_{\mu_0}$ of $\mu_0\in V^*$ given by
$
G_{\mu_0} := \{ g\in \diffG | (g^{-1})^* \mu_0 = \mu_0 \}.
$

As a consequence of the previous construction, the semi-direct product reduction theorem asserts that, for a given Hamiltonian functional $\map{H}{\gothsStar}{\bb{R}}$, the Hamiltonian dynamics on $\gothsStar = \gothgstar\times V^*$ are given by 
$$
\dot{x} = - \adSdual{\delta_xH} (x),
$$
where $\delta_xH \in \goths$ is the variational derivative of $H$ with respect to $x$, and $\adSdual{\delta_xH}$ is the dual of the adjoint operator $\adS{\delta_xH}$ of $\goths$ given in (\ref{eq:adstar_s}).

The corresponding Lie-Poisson bracket associated to $\gothsStar$ is defined by
\begin{equation}\label{eq:Lie_Posson_bracket_def}
\Pbrack{F}{G} (x) := \pair{x}{\LbrS{\varD{F}{x}}{\varD{G}{x}}}_\goths,
\end{equation}
for any $F,G \in \spFn{\gothsStar}$ and $x\in \gothsStar$ .
\end{proof}

A key advantage of the equations of motion in the Hamiltonian form (\ref{eq:Ham_kin_subs}) is that the structure of the equations clearly separates the kinetic energy functional $H_k(\alpha,\mu)$ from the underlying interconnection structure governing the evolution of the energy variables $\state$. Note that the second equation in (\ref{eq:Ham_kin_subs}) represents the advection law for the conservation of mass, whereas the term with the diamond operator in the first equation corresponds to the effect of advection on the momentum balance.

\subsection{Port-based Representation}
\newcommand{\eAlpha}{\delta_\alpha H_k}
\newcommand{\eMu}{\delta_\mu H_k}
\newcommand{\fAlpha}{\dot{\alpha}}
\newcommand{\fMu}{\dot{\mu}}
\newcommand{\eX}{(e_\alpha,e_\mu)}
\newcommand{\vecEalpha}{\hat{e}_\alpha}

In the port-based paradigm, the system described by the equations (\ref{eq:Ham_kin_subs}) is represented by two subsystems that are connected together using ports.
The first subsystem corresponds to the storage property of the system's energy (\ref{eq:reduced_Hamiltonian}), while the second subsystem corresponds to the interconnection structure encoded in the Lie-Poisson bracket (\ref{eq:Lie-Poisson_bracket_kin}).

In general, an energy-storage system (or element) in the port-Hamiltonian framework is defined by the smooth state space manifold $\cl{X}$ with a Hamiltonian functional $\map{H}{\cl{X}}{\bb{R}}$ representing the stored energy, and $x\in \cl{X}$ is called the energy variable.
The rate of change of the energy is given by
$
\dot{H} = \pair{\delta_xH}{\dot{x}}_\cl{X},
$
where $(\delta_xH,\dot{x}) \in T_x^*\cl{X}\times T_x\cl{X}$ are referred to as the effort and flow variables of the energy storage element, respectively.
The duality pairing between an effort in $T_x^*\cl{X}$ and a flow in $ T_x\cl{X}$ corresponds to the power entering the energy-storage element at a certain instant of time.

For the case of the kinetic energy Hamiltonian $H_k$ in (\ref{eq:reduced_Hamiltonian}), the state space manifold is given by $\cl{X} = \gothsStar$, the energy variables are given by $x = \state$, while the effort variable $\delta_xH_k$ and flow variable $\dot{x}$ are given by
\begin{align}
\delta_xH_k =& (\eAlpha,\eMu) \in T_{\state}^*\gothsStar \cong \gothg\times V,\\
\dot{x} =& (\fAlpha,\fMu) \in  T_{\state}\gothsStar \cong \gothgstar\times V^*.
\end{align}
The energy balance is expressed by
\begin{equation}
\begin{split}
\dot{H}_k =& \pair{\dot{x}}{\delta_xH_k}_\goths = \pair{(\fAlpha,\fMu)}{(\eAlpha,\eMu)}_\goths\\
		=& \pair{\fAlpha}{\eAlpha}_\gothg + \pair{\fMu}{\eMu}_V = \int_M \fAlpha\wedge \eAlpha + \fMu \wedge \eMu.
\end{split}
\end{equation}
Thus, this represents the first subsystem corresponding to the storage of the kinetic energy $H_k$ in (\ref{eq:reduced_Hamiltonian}).

The second subsystem corresponds to the Lie-Poisson structure defined by the map
\begin{equation}
\fullmap{J_x}{T_x^*\gothsStar \cong \goths}{T_x\gothsStar\cong \gothsStar}{e_{sk}}{J_x(e_{sk}) =: f_{sk}.}
\end{equation}
For the energy variables, $x=\state \in \gothgstar\times V^*$, we have that $e_{sk}= \eX$ and $f_{sk} = (f_\alpha,f_\mu)$. Then we can express the Lie-Poisson structure as
\begin{equation}\label{eq:Jx_1}
J_x\eX = \TwoVec{-ad^*_{e_\alpha}(\alpha) - e_\mu \diamond \mu}{-\LieD{\vecEalpha}{(\mu)}}, \qquad \eX \in \goths = \gothg\times V.
\end{equation}
By substituting the expressions of $ad^*$ in Prop. \ref{prop:adStar} and the diamond operator in Prop. \ref{prop:diamond} (case 1), we can rewrite the Lie-Poisson structure (\ref{eq:Jx_1}) as
\begin{equation}\label{eq:J_x_2}
J_x\eX = \TwoVec{-\LieD{\vecEalpha}{(\alpha)} - \divr{\vecEalpha}\alpha - (*\mu) \extd e_\mu}{- \extd((*\mu) e_\alpha)},
\end{equation}
where the second row in (\ref{eq:J_x_2}) follows from
$\LieD{\vecEalpha}{(\mu)} = \extd \iota_{\vecEalpha}(\mu) = \extd ((*\mu)\iota_{\vecEalpha}\volF) = \extd((*\mu) e_\alpha),$
using identity $\iota_{\hat{\omega}} \volF = \omega$ for any $\omega\in \gothg$, and the fact that $\mu$ is a top form.
Thus, this represents the second subsystem corresponding the interconnection structure of the system.
Since the Lie-Poisson structure (\ref{eq:J_x_2}) is defined in terms of the momentum variable $\alpha$, we refer to it as the \textbf{momentum representation} of the Lie-Poisson structure.

\begin{figure}
\centering
\includegraphics[width=0.6\textwidth]{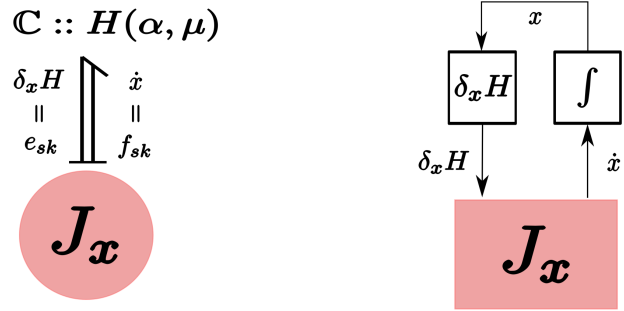}
\caption{Port-based representation of the Hamiltonian dynamics (\ref{eq:Ham_kin_subs}) corresponding to the system's kinetic energy.}
\label{fig:port-based_Ham_kinetic}
\end{figure}

The port-Hamiltonian representation of the system (\ref{eq:Ham_kin_subs}) is now given by
\begin{equation}\label{eq:port-Repr-Ham}
\dot{x} = J_x(\delta_xH_k),
\end{equation}
which is constructed by connecting the two ports $(e_{sk},f_{sk})$ and $(\delta_xH_k,\dot{x})$ together:
$$f_{sk} = \dot{x} = (\fAlpha,\fMu), \qquad e_{sk} = \delta_xH_k = (\eAlpha,\eMu).$$

Graphically, the port-based representation of the Hamiltonian dynamics (\ref{eq:Ham_kin_subs}) is shown in Fig. \ref{fig:port-based_Ham_kinetic}, where the left figure is represented using bond graphs and the right figure is represented using block diagrams, explicating causality.
The kinetic energy storage subsystem is denoted in generalized bond graphs \citep{breedveld1984physical} by a $\bb{C}$-element with its energy functional $H_k$.
The storage subsystem is connected to the Lie-Poisson structure $J_x$ through a port denoted by a double half-arrow.
The flow and effort variables are indicated on the right and left of the port, respectively.

At any instant of time, the duality pairing between the flow $\dot{x}$ and the effort $\delta_xH$ equals the power (rate of change of kinetic energy).
The conservation of energy can be seen from
\begin{align}
\dot{H}_k =& \pair{\dot{x}}{\delta_xH_k}_\goths 	= \pair{J_x(\delta_xH_k)}{\delta_xH_k}_\goths \nonumber\\
			=& \pair{-\adSdual{\delta_xH_k}(x)}{\delta_xH_k}_\goths = -\pair{x}{\adS{\delta_xH_k}{\delta_xH_k}}_\goths
												= -\pair{x}{\LbrS{\delta_xH_k}{\delta_xH_k}}_\goths = 0,	\label{eq:conserv_Ham_energy}		
\end{align}
which follows from the skew-symmetry of the Lie bracket of $\goths$.
Consequently, the Lie-Poisson structure $J_x$ is a skew-symmetric operator, which corresponds to it being a power-continuous element in port-based terminology.
It is extremely important to note that in the energy balance (\ref{eq:conserv_Ham_energy}), the surface terms (\ref{eq:pairing_adSdual}) that should appear in the fourth equality naturally disappear for the closed system (\ref{eq:Ham_kin_subs}).


In summary, the closed port-Hamiltonian system  shown in Fig. \ref{fig:port-based_Ham_kinetic} describes the conservation of kinetic energy $H_k$ and the corresponding evolution of the energy variables $\state$.
The conservation of energy follows from the skew-symmetry of the Lie-Poisson structure.
The aforementioned port-Hamiltonian system is still equivalent to the \textit{standard} Hamiltonian one that describes a conservative system that is isolated from any energy-exchange with the external world.
Next, we discuss how to allow non-zero energy exchange by replacing the underlying Lie-Poisson structure with a \textit{Stokes-Dirac structure}.

\subsection{Open Model of Kinetic Energy}
\newcommand{\fstr}{\text{f}_\text{s}}
\newcommand{\eBoud}{e_{\partial k}}
\newcommand{\fBoud}{f_{\partial k}}

There are two ways in which the port-Hamiltonian system (\ref{eq:port-Repr-Ham}) can interact and exchange energy with the world; either through the boundary $\Mbound$ of the spatial manifold $M$ or within the domain itself through a distributed port that allows energy exchange at every point in $M$.
The former allows exchange of kinetic energy by mass inflow or outflow, while the latter allows transformation of kinetic energy to another form in a reversible (or irreversible) way.

To add a distributed port based on Newton's second law, a \textit{distributed force} field $\fstr \in \gothgstar = \spKFormM{1}$ is added to the momentum balance equation such that the port-Hamiltonian system (\ref{eq:Ham_kin_subs}) is rewritten, using (\ref{eq:J_x_2}) as
\begin{align}
\TwoVec{\fAlpha}{\fMu} =&  \TwoVec{-\LieD{(\eAlpha)^\wedge}{(\alpha)} - \divr{(\eAlpha)^\wedge}\alpha - (*\mu) \extd \eMu + \fstr}{- \extd((*\mu) \eAlpha)},\label{eq:pH_kinetic_sys}\\
H_k\state =& \int_M \frac{1}{2(*\mu)} \alpha \wedge * \alpha.\label{eq:pH_kinetic_Ham}
\end{align}
It is worth noticing how our previous choice of the Lie algebra $\gothg$ as $\spKFormM{n-1}$ effects the external force field $\fstr$ on the dual algebra $\gothg^*$ to be a co-vector field, which is the correct geometrical representation of a force field.
We can write (\ref{eq:pH_kinetic_sys}) more compactly as
\begin{equation}\label{eq:pH_kinetic_sys_compact}
\dot{x} = J_x(\delta_x H_k) + G \fstr,
\end{equation}
with $J_x$ given by (\ref{eq:J_x_2}) and $\map{G}{\gothgstar}{\gothgstar\times V^*}$, $G=(1 \quad 0)^\top$ representing the input map.
The force one-form (co-vector field) $\fstr$ will be used later in Part II to model stress forces due to pressure.
In general, the distributed force can be used for modeling other stress forces due to viscosity, as well as any body (external) forces on the continuum (e.g. due to magnetic fields or gravity and electrostatic accelerations).

The Hamiltonian energy (\ref{eq:pH_kinetic_Ham}), as a functional $\map{H_k}{\gothsStar}{\bb{R}}$, admits its rate of change such that along trajectories $x(t)$, parameterized by time $t\in \bb{R}$, it holds that
\begin{equation}\label{eq:energy_balance_Hk}
\dot{H}_k=  \pair{\dot{x}}{\delta_xH_k}_\goths = \int_{M} \dot{x} \wedge \delta_{x}H_k.
\end{equation}
As shown in (\ref{eq:conserv_Ham_energy}), for an isolated fluid system on a closed manifold (corresponding to $G=0$, and either $\partial M=\emptyset$ or $\partial M$ is impermeable) the kinetic energy is always conserved.
However, for a general open fluid system, the expression for the kinetic energy balance (\ref{eq:energy_balance_Hk}) is given by the following result.

\begin{theorem}\label{theorem:energy_balance_kinetic}
The rate of change of the Hamiltonian (\ref{eq:pH_kinetic_Ham})  along trajectories of the port-Hamiltonian system (\ref{eq:pH_kinetic_sys_compact}) is given by
\begin{equation}
\label{eq:power_balance_kinetic_system}
\dot{H}_k=  \int_{\Mbound} \eBoud \wedge \fBoud   + \int_M e_d \wedge f_d,
\end{equation}
where the boundary port variables $\eBoud,\fBoud \in \spKForm{0}{\Mbound} \times \spKForm{n-1}{\Mbound}$ and distributed port variables $e_d,f_d \in \spKFormM{1} \times \spKFormM{n-1}$ are defined by
\begin{align*}
\eBoud &:=  {\bound{\left(\frac{\iota_{(\eAlpha)^\wedge}(\alpha)}{*\mu} + \eMu \right )}}  ,\quad &  e_d &:= \normalfont{\fstr},\\
\fBoud &:= {- \bound{(*\mu)\eAlpha}},\quad & f_d &:=  \eAlpha,
\end{align*}
where $\bound{\eta} \in \spKForm{k}{\Mbound} $ denotes the trace of the form $\eta \in \spKFormM{k}$. 
\end{theorem}
\begin{proof}
For notational simplicity, we denote $\eX = (\eAlpha,\eMu)$.
By substituting (\ref{eq:pH_kinetic_sys_compact}) in (\ref{eq:energy_balance_Hk}), we have that
\begin{align*}
\dot{H}_k 	=& \pair{J_x\eX}{\eX}_\goths + \pair{G \fstr}{\eX}_\goths\\
			=& -\pair{\adSdual{\eX}\state}{\eX}_\goths + \pair{\fstr}{G^\top\eX}_\goths\\ 
			=& -\underbrace{\pair{\state}{\adS{\eX}\eX}_\goths}_{=0} -  \int_{\Mbound}\bound{ \surfAdSstar{\eX}(\alpha,\mu,e_\alpha,e_\mu)} + \pair{\fstr}{e_\alpha}_\gothg, 																																				
\end{align*}
where (\ref{eq:pairing_adSdual}) was used, and  the first term in the last equality vanishes due to the skew-symmetry property of the Lie bracket $\adS{} = \LbrS{}{}$.
Using (\ref{eq:surfAdS}),(\ref{eq:surfAd}), and Prop. \ref{prop:diamond} (case i),  we can express the surface term as
{
\begin{align*}
\surfAdSstar{\eX}(\alpha,\mu,e_\alpha,e_\mu) 	&= \surfAdStar{e_\alpha}(\alpha,e_\alpha) -\surfPhiDTA{\mu}(e_\alpha,e_\mu) - \surfLie{\hat{e}_\alpha}{\mu}{e_\mu}  + \surfLie{\hat{e}_\alpha}{\mu}{e_\mu} ,\\
			&= \iota_{\hat{e}_\alpha}(\alpha\wedge e_\alpha) - (-(*\mu)e_\alpha\wedge e_\mu ) + 0.
\end{align*}
Using the interior product properties, we have that
\begin{align*}
\iota_{\hat{e}_\alpha}(\alpha\wedge e_\alpha) &= \iota_{\hat{e}_\alpha}(\alpha)\wedge e_\alpha - \alpha\wedge \iota_{\hat{e}_\alpha}(e_\alpha) \\
			&= \iota_{\hat{e}_\alpha}(\alpha)\wedge e_\alpha - \alpha\wedge \iota_{\hat{e}_\alpha}\circ\iota_{\hat{e}_\alpha}(\volF) 
			=  \iota_{\hat{e}_\alpha}(\alpha)\wedge e_\alpha,						
\end{align*}
where the nil-potency property of the interior product was used in the last equality.
}
Therefore, finally we have that
\begin{equation}\label{eq:energy_balance_old}
\dot{H}_k = \int_{\Mbound} { \bound{\left(\frac{\iota_{\hat{e}_\alpha}(\alpha)}{*\mu} + e_\mu \right)}  \wedge  \bound{-(*\mu)e_\alpha}}+ \int_{M} \fstr \wedge e_\alpha,
\end{equation}
which concludes the proof.
\end{proof}

It is interesting to note the physical interpretation of the boundary { flow variable $\fBoud$}. Using (\ref{eq:var_der_H_k}) we have that
$$\fBoud = - \bound{(*\mu)\eAlpha} = - \bound{*\alpha} = -\bound{(*\mu) \omega_v}= -\bound{\iota_v\mu},$$
which represents the incoming mass flow through the boundary. As discussed before in Sec. \ref{sec:permeable_vs_imperm_boundaries}, the condition $\bound{\omega_v} = 0$ corresponds to an impermeable boundary $\Mbound$ which implies that $\fBoud=0$ and thus results in no exchange of power through the boundary, as seen in (\ref{eq:power_balance_kinetic_system}).
The physical interpretation of the boundary effort variable $\eBoud$ will be discussed later.

\begin{figure}
\centering
\begin{subfigure}{0.45\textwidth}
\includegraphics[height=0.6\textwidth,width=\textwidth]{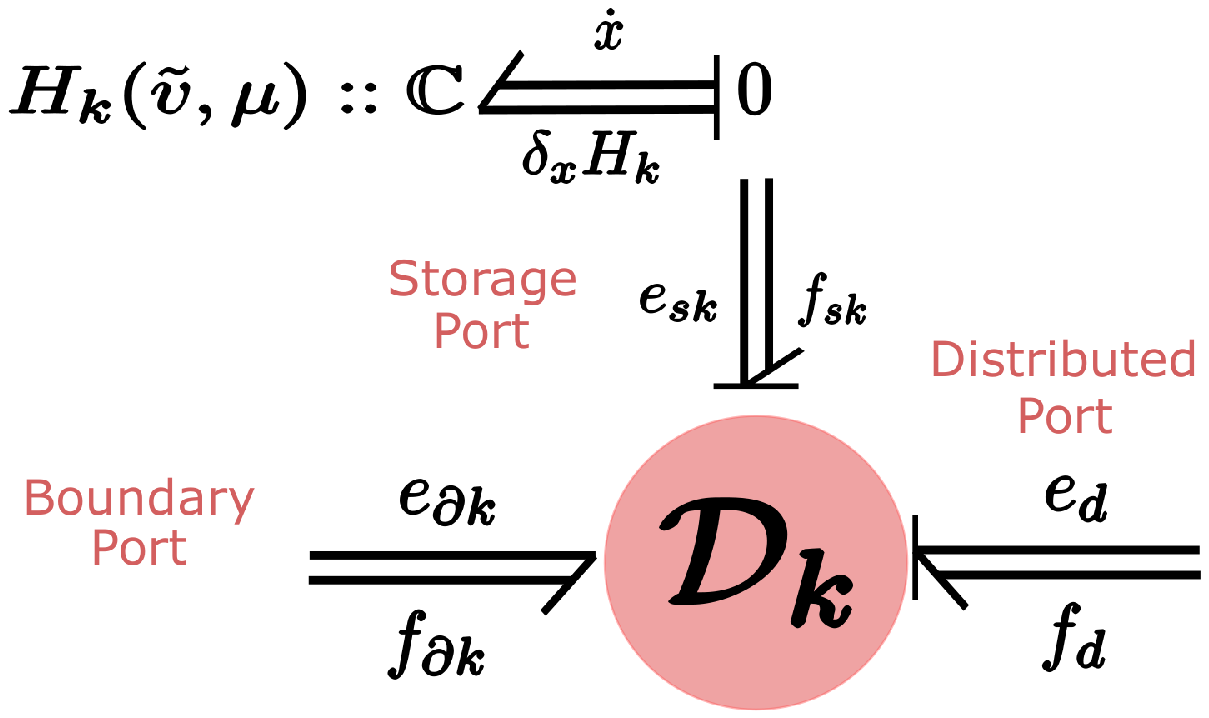}
\caption{Bond Graph}
\end{subfigure}
\begin{subfigure}{0.45\textwidth}
\includegraphics[height=0.7\textwidth,width=\textwidth]{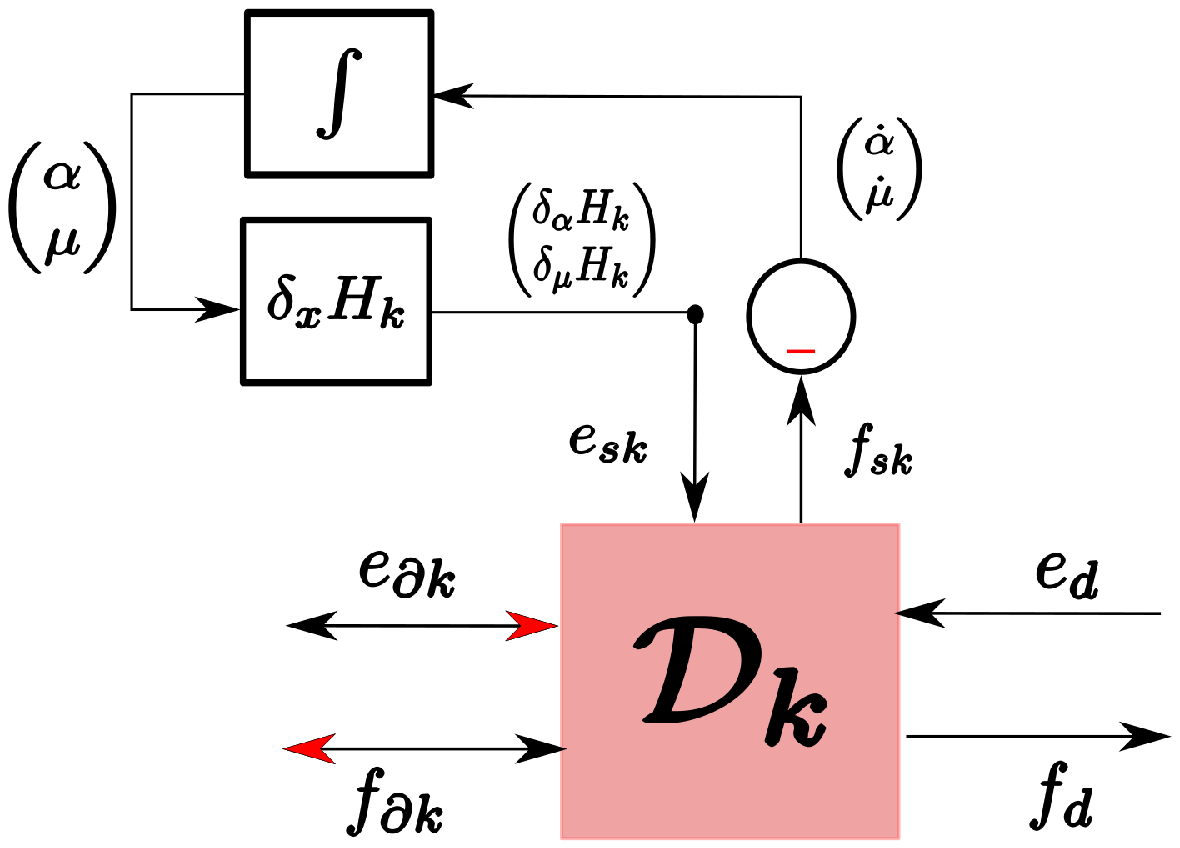}
\caption{Block Diagram}
\end{subfigure}
\caption{Graphical Representation of the Kinetic Energy Subsystem (\ref{eq:pH_kinetic_sys}) in terms of the state variables: momentum one-form $\alpha$ and the mass top-form $\mu$.}
\label{fig:SDS_Kinetic_Energy_Mom}
\end{figure}

With reference to Fig. \ref{fig:SDS_Kinetic_Energy_Mom}, the port-Hamiltonian system now has 3 ports; an energy storage port $(e_{sk},f_{sk})$, a boundary port $(\eBoud,\fBoud)$, and a distributed port $(e_d,f_d)$.
The port $(\eBoud,\fBoud)$ represents the power incoming to the system through the boundary, while the port $(e_d,f_d)$ represents the power flow from an external distributed source in the spatial domain.
The 3 ports are connected together through a power-conserving structure $\cl{D}_k$ that encodes the power balance equation (\ref{eq:power_balance_kinetic_system}).
Such mathematical structure is called a \textbf{Stokes-Dirac structure} defined as
\begin{equation}\label{eq:SDS_kinetic_system}
\begin{split}
\cl{D}_k = \{ (f_{sk},&\fBoud,f_d,e_{sk},\eBoud,e_d) \in \cl{B}_k | \\
			 \TwoVec{f_\alpha}{f_\mu} &= \TwoVec{\LieD{\vecEalpha}{(\alpha)} + \divr{\vecEalpha}\alpha + (*\mu) \extd e_\mu}{ \extd((*\mu) e_\alpha)} - \TwoVec{1}{0 }e_d,\\
			f_d&= \begin{pmatrix} 1 &  0\end{pmatrix} \TwoVec{e_\alpha}{e_\mu},\\ 
			 \TwoVec{\eBoud}{\fBoud}  &= \TwoVec{\bound{\left(\frac{\iota_{\hat{e}_\alpha}(\alpha)}{*\mu} + e_\mu \right)}}{-(*\mu)e_\alpha\bound{}}\},
\end{split}
\end{equation}
where the storage port variables are given by $f_{sk} = ({f_\alpha},{f_\mu}) \in \gothsStar = \spKFormM{1} \times \spKFormM{n} $ and $e_{sk} = (e_{\alpha},e_{\mu}) \in \goths = \spKFormM{n-1} \times \spKFormM{0} $.
The bond-space $\cl{B}_k= \cl{F}_k\times \cl{E}_k$ is the product space of the flow space $\cl{F}_k=  \spKFormM{1} \times \spKFormM{n} \times \spKForm{n-1}{\partial M} \times \spKForm{n-1}{M} $ and the effort space $\cl{E}_k = \spKFormM{n-1} \times \spKFormM{0}  \times \spKForm{0}{\partial M} \times \spKForm{1}{M}$, represented by smooth differential forms of the appropriate degree.

The Dirac structure (\ref{eq:SDS_kinetic_system}) is a modulated one, in the sense that it depends on the state variables $x=\state$.
In the absence of the distributed force ($e_d = 0$) and in case of an impermeable (or no) boundary ($\fBoud = 0$), the Stokes-Dirac structure degenerates to the Lie-Poisson structure (\ref{eq:J_x_2}) which encodes the conservation of kinetic energy ($\dot{H}_k=0$).
The sign difference between the Lie-Poisson structure in (\ref{eq:J_x_2}) and its counterpart in (\ref{eq:SDS_kinetic_system}) is due to the choice of having all power ports entering the Dirac structure.
This can be clearly seen in the graphical representation of the port-Hamiltonian system in Fig. \ref{fig:SDS_Kinetic_Energy_Mom} in terms of bond graph and block diagram notation.

It is interesting to note that the three power ports shown in Fig. \ref{fig:SDS_Kinetic_Energy_Mom} are of different natures and result from different duality pairings.
The power in the storage port is expressed by the pairing $\pair{e_{sk}}{f_{sk}}_\goths$, the power in the distributed port is expressed by the pairing $\pair{e_d}{f_d}_\gothg$, while the power in the distributed port is expressed by $\int_{\Mbound} \eBoud\wedge\fBoud$.


The port-Hamiltonian system (\ref{eq:pH_kinetic_sys}) can be recovered from the Stokes-Dirac structure (\ref{eq:SDS_kinetic_system}) by imposing 
\begin{align*}
e_{sk} &= \TwoVec{e_{\alpha}}{e_{\mu}}=\TwoVec{\eAlpha}{\eMu}, \qquad & f_{sk} =& \TwoVec{f_{\alpha}}{f_{\mu}}=\TwoVec{-\fAlpha}{-\fMu},\\
e_d &= \fstr ,\qquad & f_d =& \eAlpha,
\end{align*}
where the minus sign is due to the choice of having the storage port entering the Dirac structure, as shown in Fig. \ref{fig:SDS_Kinetic_Energy_Mom}.
Therefore the implicit port-Hamiltonian dynamics are governed by
\begin{equation}
((-\fAlpha,-\fMu),\fBoud,\eAlpha,(\eAlpha,\eMu),\eBoud,\fstr) \in \cl{D}_k.
\end{equation}

\subsection{Change of Coordinates to the Velocity Representation}
\newcommand{\stateV}{(\vf,\tilde{\mu})}
\newcommand{\xT}{\tilde{x}}
\newcommand{\eVt}{\delta_{\vf} H_k}
\newcommand{\eMuT}{\delta_{\tilde{\mu}} H_k}
\newcommand{\fVt}{\dot{\vf}}
\newcommand{\fMuT}{\dot{\tilde{\mu}}}
\newcommand{\eXt}{(e_{\vf},e_{\tilde{\mu}})}
\newcommand{\JxT}{\tilde{J}_{\xT}}
\newcommand{\muT}{{\tilde{\mu}}}

\newcommand{\eSV}{{e}_{sk}}
\newcommand{\fSV}{{f}_{sk}}

In the work of \cite{van2002hamiltonian}, the port-Hamiltonian model for compressible isentropic flow was given in terms of the velocity 1-form $\vf \in \spKFormM{1}$ and the mass top-form $\mu \in \spKFormM{n}$ as the underlying energy variables.
The arguments behind this choice of coordinates where not connected to the geometric structure underlying the state space, as we presented so far in this work.
Instead, in \citep{van2002hamiltonian} the authors defined a canonical version of a Stokes Dirac structure based on exterior derivatives, that is representative of a number of physical systems based on conservation laws of different nature.
Then, the authors modified this canonical Dirac structure to represent the fluid dynamical system with the velocity and mass forms as energy variables.

One of the contributions of this article is to present a rigorous derivation of the Stokes Dirac structure presented in \citep{van2002hamiltonian}, starting from (\ref{eq:SDS_kinetic_system}) based on semi-direct product theory.
We will represent the port-Hamiltonian system in (\ref{eq:pH_kinetic_sys}) in terms of the same energy variables in \citep{van2002hamiltonian} by means of a change of coordinates, recovering the Dirac structure which was introduced there without a formal derivation.

The following change of coordinates procedure was introduced by \cite{Vankerschaver2010} for the Lie-Poisson structure only.
In this work we extend the coordinate change for the full Stokes-Dirac structure and emphasize its port-based interpretation.
It will turn out that in terms of the new variables, the expression of the Dirac structure will be simpler than the momentum representation in (\ref{eq:J_x_2}). In addition, this description yields further advantages such as a trivial derivation of the vorticity equation.

Consider the (nonlinear) diffeomorphism 
\begin{equation}\label{eq:Vel_diffeom}
\Phi:\state\mapsto\stateV:= (\frac{\alpha}{*\mu},\mu).
\end{equation}
We denote the new state variables $\stateV$ by $\xT \in \gothsStar$, i.e. $\Phi:x\mapsto\xT$.
Consider the pushforward map $(\Phi_*)_x$ and pullback map $(\Phi^*)_x$ of $\Phi$ at a point $x\in \gothsStar$
\begin{equation}
\fullmap{(\Phi_*)_x}{T_x\gothsStar}{T_{\xT}\gothsStar}{ (\fAlpha,\fMu)}{ (\fVt,\fMuT),}
\qquad \qquad 
\fullmap{(\Phi^*)_x}{T_{\xT}^*\gothsStar}{T_x^*\gothsStar}{\eXt}{ \eX,}
\end{equation}



Now we consider the coordinate change for the Lie-Poisson part (\ref{eq:J_x_2}) of the Stokes-Dirac structure.
The \textbf{velocity representation} of the Lie-Poisson structure is the map
$
\map{\JxT}{T_{\xT}^*\gothsStar}{T_{\xT}\gothsStar},
$
defined by
\begin{equation}\label{eq:JxT_def}
\JxT := (\Phi_*)_x\circ J_x \circ (\Phi^*)_x.
\end{equation}
The change of coordinates described above can be represented graphically in an elegant way in the port-based framework as shown in the Fig. \ref{fig:coordinates_LiePoisson}.

\begin{figure}
\centering
\includegraphics[width=\textwidth]{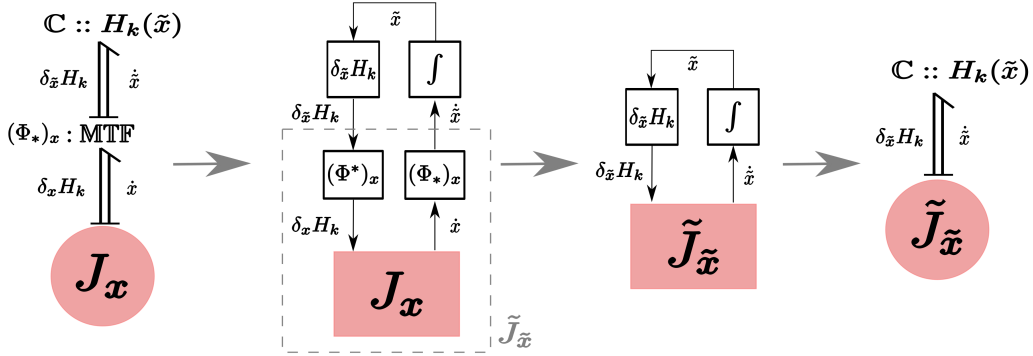}
\caption{Graphical representation of change of coordinates from the momentum-representation $J_x$ to the velocity representation $\JxT$ of the Lie-Poisson structure.}
\label{fig:coordinates_LiePoisson}
\end{figure}

The exact expressions for the maps $(\Phi_*)_x, (\Phi^*)_x,$ and $\JxT$ are given by the following theorem \citep{Vankerschaver2010}.
\begin{proposition}\label{proposition:vel_represent_maps}
The pushforward $(\Phi_*)_x$ and pullback map $(\Phi^*)_x$ of the diffeomorphism (\ref{eq:Vel_diffeom}) at a point $x = \state \in \gothsStar$ are given respectively by
\begin{align}
(\Phi_*)_x(\fAlpha,\fMu) &= \left( \frac{\fAlpha - (*\fMu) \vf}{*\muT},\fMu \right) =: (\fVt,\fMuT),\label{eq:pushfwd_vel_map}\\
(\Phi^*)_x\eXt &= \left( \frac{e_{\vf}}{*\muT}, e_{\muT} - \frac{*(\vf\wedge e_{\vf})}{*\muT} \right) =: \eX,\label{eq:pullback_vel_map}
\end{align}
where $\stateV = \Phi\state$.
The velocity representation of the Lie-Poisson structure defined by (\ref{eq:JxT_def}) is expressed as
\begin{equation}\label{eq:Jx_tilde}
\JxT\eXt = \TwoVec{-\extd e_{\muT} - \frac{1}{*\muT} \iota_{\hat{e}_{\vf}} \extd \vf }{-\extd e_{\vf} }.
\end{equation}
\end{proposition}
\begin{proof}
See Appendix Sec. \ref{proof:vel_represent_maps}.
\end{proof}

\renewcommand{\stateV}{(\vf,{\mu})}
In terms of the new energy variables, the kinetic energy Hamiltonian functional and its variational derivatives are given by the following result.
\begin{proposition}\label{prop:kinetic_Ham_new}
The kinetic energy Hamiltonian given by (\ref{eq:reduced_Hamiltonian}) in terms of the velocity 1-form $\vf$ and the mass top-form $\mu$ is given by
\begin{equation}\label{eq:new_Hamiltonian}
H_k(\vf,\mu) = \int_M \half (*\mu) \vf \wedge * \vf.
\end{equation}
The variational derivatives $\delta_{\vf} H_k \in T_{\vf}^*\gothgstar \cong \gothg = \spKFormM{n-1}$ and $\delta_{\mu} H_k \in T_{\mu}^*V^* \cong V =  \spKFormM{0}$ with respect to the states $\vf \in \gothgstar $ and $\mu \in V^*$, respectively, are given by
\begin{equation}\label{eq:varDeriv_stateV}
\delta_{\vf} H_k  = (*\mu) * \vf = \iota_v \mu, \qquad\qquad \delta_{\mu} H_k = \half \iota_v \vf,
\end{equation}
\end{proposition}
\begin{proof}
See Appendix Sec. \ref{{proof:kinetic_Ham_new}}.
\end{proof}

The new energy balance for the kinetic energy Hamiltonian (\ref{eq:new_Hamiltonian}) instead of (\ref{eq:power_balance_kinetic_system}) can be now calculated using the pullback map $(\Phi^*)_x$ in (\ref{eq:pullback_vel_map}) as follows.
\begin{theorem}\label{theorem:new_power_balance_H_k}
The rate of change of the kinetic energy Hamiltonian (\ref{eq:new_Hamiltonian}) in terms of the new state variables $\stateV$ is given by
\begin{equation}
\label{eq:new_power_balance}
\dot{H}_k=  \int_{\Mbound} \eBoud \wedge \fBoud   + \int_M e_d \wedge f_d,
\end{equation}
where the boundary port variables $\eBoud,\fBoud \in \spKForm{0}{\Mbound} \times \spKForm{n-1}{\Mbound}$ and the distributed port variables $(e_d,f_d) \in \spKFormM{1} \times \spKFormM{n-1}$ are expressed using the new coordinates $\stateV$ by
\begin{align*}
\eBoud &=  { \bound{\eMu} = \bound{\half \iota_v\vf}},\quad &  e_d &= \normalfont{\fstr},\\
\fBoud &=  { - \bound{\eVt} = -\bound{\iota_v\mu}} ,\quad & f_d &=  \frac{\eVt}{*\mu} = \omega_v,
\end{align*}
where the variational derivatives of $H_k$ in (\ref{eq:new_Hamiltonian}) are given by (\ref{eq:varDeriv_stateV}), and $\omega_v = \iota_v\volF \in \spKFormM{n-1}$.
\end{theorem}
\begin{proof}
See Appendix Sec. \ref{proof:new_power_balance_H_k}
\end{proof}

As mentioned before, the boundary flow variable $\fBoud$ is physically the incoming mass flow through the boundary. On the other hand, the boundary effort variable $\eBoud$ is physically the dynamic pressure (modulo the density) which is the only pressure that exists related to the kinetic energy of the fluid since we have not modeled yet any thermodynamic potential.

\newcommand{\DkT}{\cl{\tilde{D}}_k}
The new Stokes-Dirac structure $\DkT$ in terms of the state variables $\stateV$ can now be constructed using the velocity representation of the Lie-Poisson structure (\ref{eq:Jx_tilde}) and the definitions of the boundary and distributed port variables in Theorem \ref{theorem:new_power_balance_H_k}.
The expression for $\DkT$ is given by
\begin{equation}\label{eq:SDS_kinetic_system_new}
\begin{split}
\DkT = \{ (\fSV, &\fBoud, f_d,\eSV,\eBoud,e_d) \in \cl{B}_k | \\
			 \TwoVec{f_\vf}{f_\mu} &= \TwoVec{\extd e_\mu +  \frac{1}{*\mu} \iota_{\hat{e}_{\vf}} \extd \vf }{\extd e_\vf }- \TwoVec{\frac{1}{*\mu}}{0 }e_d,\\
			f_d&= \begin{pmatrix} \frac{1}{*\mu} &  0\end{pmatrix} \TwoVec{e_\vf}{e_\mu},\\ 
			 \TwoVec{\eBoud}{\fBoud}  &=  \TwoTwoMat{0}{1}{-1}{0}\TwoVec{\bound{e_\vf}}{\bound{e_\mu}}\},
\end{split}
\end{equation}
 where the storage port variables are now given by $f_{sk} = (f_\vf,f_\mu) \in \gothsStar = \spKFormM{1} \times \spKFormM{n} $ and $e_{sk} = (e_\vf,e_\mu) \in \goths = \spKFormM{n-1} \times \spKFormM{0} $.
The bond-space $\cl{B}_k$ is the same as the one defined before for $\cl{D}_k$ in (\ref{eq:SDS_kinetic_system}).

By excluding the distributed port $(e_d,f_d)$, the Stokes-Dirac structure (\ref{eq:SDS_kinetic_system_new}) coincides with the one postulated in \citep{van2002hamiltonian} for an ideal fluid and considered as an extended version of with respect to the one containing only exterior derivative operators.


\begin{figure}
\centering
\begin{subfigure}{0.45\textwidth}
\includegraphics[height=0.6\textwidth,width=\textwidth]{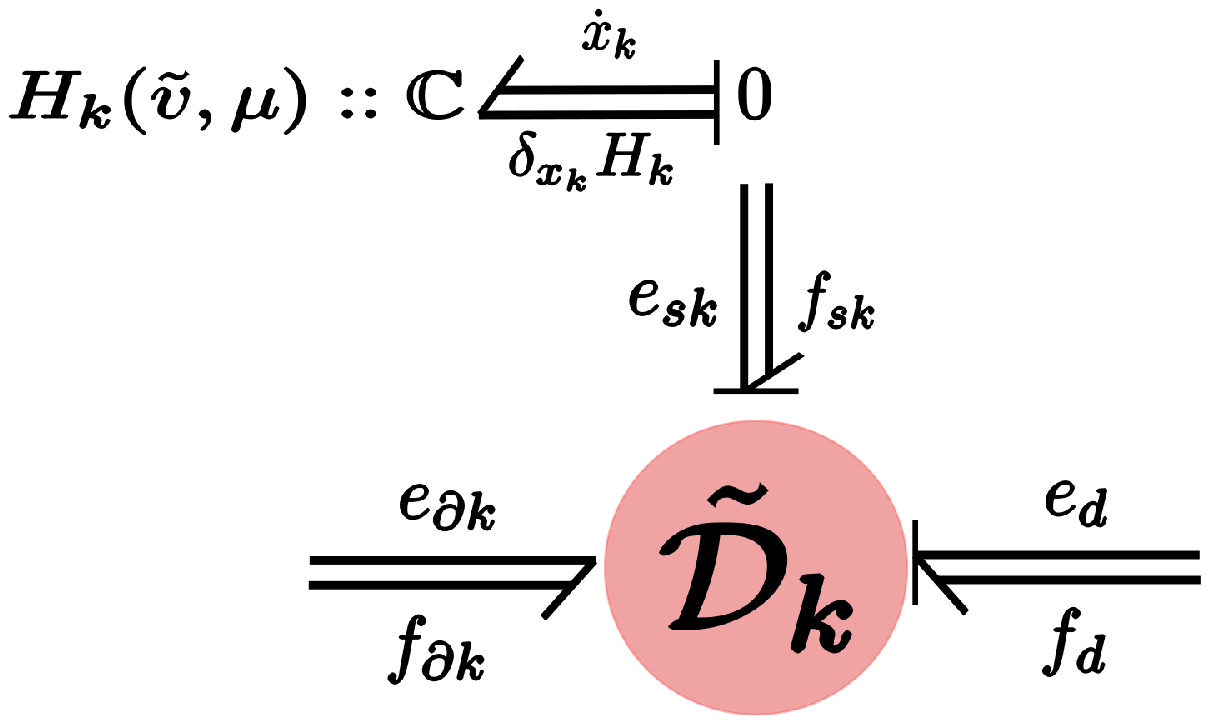}
\caption{Bond Graph}
\end{subfigure}
\begin{subfigure}{0.45\textwidth}
\includegraphics[height=0.7\textwidth,width=\textwidth]{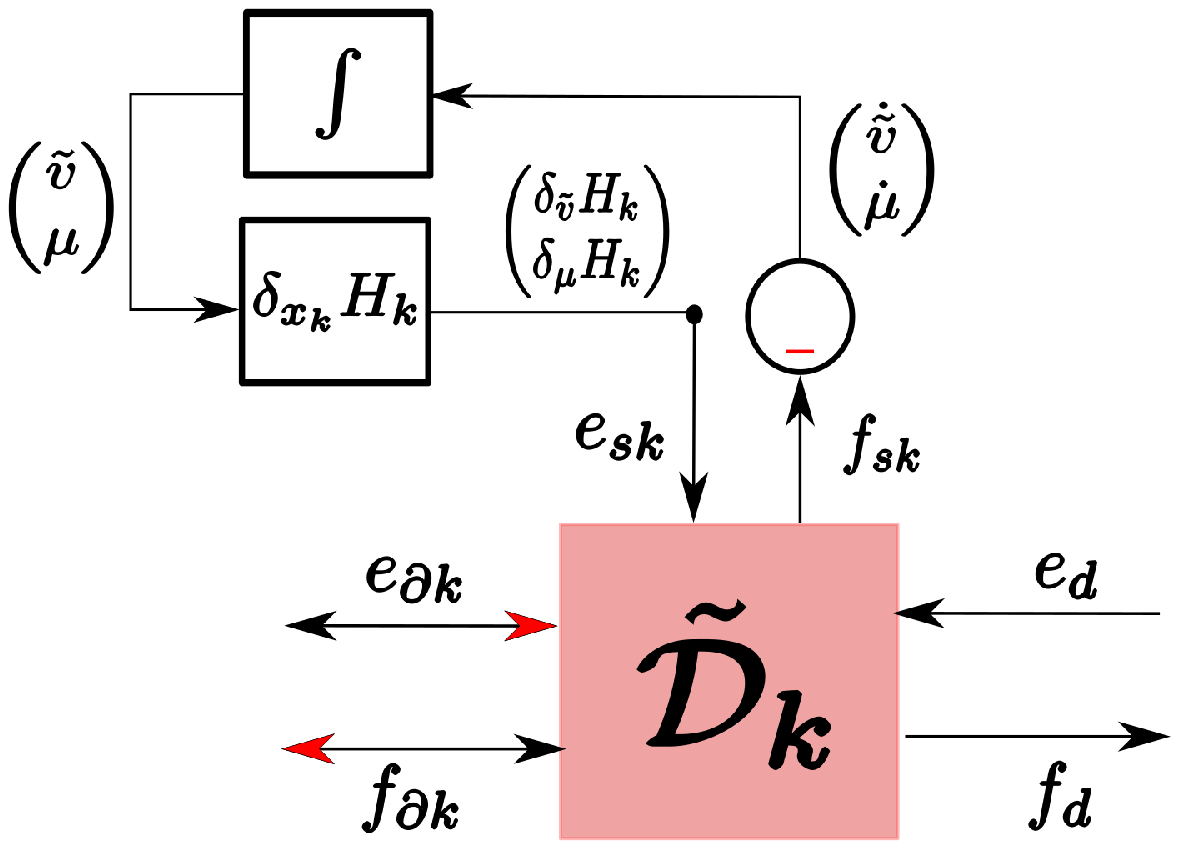}
\caption{Block Diagram}
\end{subfigure}
\caption{Graphical Representation of the Kinetic Energy Subsystem (\ref{eq:pH_kinetic_sys_new_1}) in terms of the state variables: velocity one-form $\vf$ and the mass top-form $\mu$.}
\label{fig:SDS_Kinetic_Energy_Vel}
\end{figure}

To summarize, the explicit port-Hamiltonian dynamical model in terms of the kinetic energy state variable $x_k:= \xT = \stateV \in \cl{X} = \gothsStar$ is given by
\begin{align}
\TwoVec{\fVt}{\fMu} =& \TwoVec{-\extd (\eMu) -  \iota_{v} \extd \vf } {-\extd (\eVt)} + \TwoVec{ \frac{1}{*\mu}}{0}\fstr, \label{eq:pH_kinetic_sys_new_1}\\
\omega_v =& \begin{pmatrix}  \frac{1}{*\mu} &  0\end{pmatrix} \TwoVec{\eVt}{\eMu}, \label{eq:pH_kinetic_sys_new_2}\\
H_k(x_k) =& H_k\stateV = \int_M \half (*\mu) \vf \wedge * \vf,\label{eq:pH_kinetic_Ham_new}
\end{align}
where the variational derivatives are given by Prop. \ref{prop:kinetic_Ham_new}.
The term $\iota_{v} \extd \vf $ is equivalent to $\frac{1}{*\mu} \iota_{\hat{e}_{\vf}} \extd \vf = \iota_{({ \hat{e}_{\vf}}/{*\mu} )}\extd \vf$ in (\ref{eq:Jx_tilde}), which follows from
$
\frac{ \hat{e}_{\vf}}{*\mu} = \frac{(\eVt)^\wedge}{*\mu} = \frac{(*\mu *\vf)^\wedge}{*\mu} = (*\vf)^\wedge = v.
$
Moreover, the port-Hamiltonian dynamics (\ref{eq:pH_kinetic_sys_new_1}-\ref{eq:pH_kinetic_Ham_new}) are recovered from the Dirac structure $\DkT$ in (\ref{eq:SDS_kinetic_system_new}) by
\begin{equation}
((-\fVt,-\fMu),\fBoud,\omega_v,(\eVt,\eMu),\eBoud,\fstr) \in \DkT.
\end{equation}
The graphical representation of the kinetic energy port-Hamiltonian system (\ref{eq:pH_kinetic_sys_new_1}-\ref{eq:pH_kinetic_Ham_new}) is shown in Fig. \ref{fig:SDS_Kinetic_Energy_Vel}.

%

\section{Conclusion}

In part I of this paper, we presented the systematic port-Hamiltonian formulation of the kinetic energy storage subsystem and the derivation of its corresponding Stokes-Dirac structure that routes the energy exchanged with a distributed port (which allows the connection to subsystems that interact with the 
fluid at every point of its domain) and a boundary port (which allows the connection to subsystems that allow exchange of mass flow).

Using the geometric nature and the semi-direct product algebra structure of the state space underlying the energy variables, we showed that the Stokes-Dirac structure could be systematically constructed in two steps.
First one uses Poisson reduction \citep{marsden1984semidirect,marsden1984reduction} to construct the underlying Lie-Poisson structure and its corresponding bracket.
The second step then extends the Lie-Poisson structure to a Stokes-Dirac structure that allows non-zero exchange of kinetic energy within the spatial domain and through its boundary.
All the surface terms were easily identified thanks to the duality pairings of the different maps constructed on the semi-direct product algebra $\goths$.

This systematic two-step procedure coincide with the Stokes-Dirac structure of compressible isentropic flow heuristically obtained in \citep{van2001fluid,van2002hamiltonian} based on the velocity of the flow as the energy variable.
We also presented an alternative Stokes-Dirac structure based on the momentum of the flow, and showed how to relate both representations to each other.

The constructions of this paper, which only comprise what is needed to describe an universal subsystem for the kinetic energy storage of a fluid and its corresponding Dirac structure, will be combined in part II of this series with corresponding constructions for the description of a universal potential energy subsystem, by whose addition one obtains a number of Euler equations.

\section*{Funding}
This work was supported by the PortWings project funded by the European Research Council [Grant Agreement No. 787675]

\renewcommand\thesection{A}
\section{Appendix: Proof of Theorems and Propositions}
\subsection{Proof of Theorem \ref{theorem:adStarS_duality}}\label{proof:adStarS_duality}

Consider $\elS{1}, \elS{2} \in \goths$. Using the definition of the bracket in (\ref{eq:Liebracket_s}), the duality pairing in (\ref{eq:duality_s}), and the identity
$$\int_M a \wedge \LieD{\omV}{\abar} + \LieD{\omV}{a} \wedge \abar = \int_{\Mbound} \surfLie{\omV}{a}{\abar},$$
we have that
\begin{align*}
\pair{\elSdual{}}{\adS{\elS{1}}\elS{2}}_\goths		&= \pair{\elSdual{}}{(\LbrG{\omega_1}{\omega_2},\LieD{\hat{\omega}_2}{\abar_1}-\LieD{\hat{\omega}_1}{\abar_2}}_\goths \\
													&= \pair{\alpha}{ad_{\omega_1}(\omega_2)}_\gothg + \pair{a}{\LieD{\hat{\omega}_2}{\abar_1}}_V - \pair{a}{\LieD{\hat{\omega}_1}{\abar_2}}_V\\
													&= \pair{\alpha}{ad_{\omega_1}(\omega_2)}_\gothg - \pair{\LieD{\hat{\omega}_2}{a}}{\abar_1}_V + \pair{\LieD{\hat{\omega}_1}{a}}{\abar_2}_V\\
													&\quad + \int_{\Mbound} \surfLie{\hat{\omega}_2}{a}{\abar_1}\bound{}  - \int_{\Mbound} \surfLie{\hat{\omega}_1}{a}{\abar_2}\bound{}.
\end{align*}

Let $c := k(n-k)\in \bb{R}$. Using the wedge product properties and (\ref{eq:pair_Rmap}), the term $-\pair{\LieD{\hat{\omega}_2}{a}}{\abar_1}_V$ can be expressed as
\begin{align}
-\int_M \LieD{\hat{\omega}_2}{a} \wedge \abar_1 &= \int_M (-1)^{c+1} \abar_1 \wedge \LieD{\hat{\omega}_2}{a}  = (-1)^{c+1} \pair{\abar_1}{\mapR(\omega_2)}_{V^*} \nonumber\\
& = (-1)^{c+1} \left(\pair{\mapRdual(\abar_1)}{\omega_2}_\gothg -\int_{\Mbound} \surfPhiDTA{a}(\omega_2,\abar_1)\bound{}\right )\nonumber\\
& = \pair{\abar_1\diamond a}{\omega_2}_\gothg - \int_{\Mbound} (-1)^{c+1}  \surfPhiDTA{a}(\omega_2,\abar_1)\bound{}. \label{eq:proof_adSt_1}
\end{align}

Using the pairing equalities (\ref{eq:adStar_pairing}) and (\ref{eq:proof_adSt_1}), we have that
\begin{align*}
\pair{\elSdual{}}{\adS{\elS{1}}\elS{2}}_\goths& = \pair{ad^*_{\omega_1}(\alpha)}{\omega_2}_\gothg - \int_{\Mbound}\surfAdStar{\omega_1}(\alpha,\omega_2)\bound{} + \pair{\abar_1\diamond a}{\omega_2}_\gothg \\
& - \int_{\Mbound} (-1)^{c+1}  \surfPhiDTA{a}(\omega_2,\abar_1)\bound{} + \pair{\LieD{\hat{\omega}_1}{a}}{\abar_2}_V \\
& + \int_{\Mbound} \surfLie{\hat{\omega}_2}{a}{\abar_1}\bound{}  - \int_{\Mbound} \surfLie{\hat{\omega}_1}{a}{\abar_2}\bound{}\\
& = \pair{ad^*_{\omega_1}(\alpha) + \abar_1 \diamond a}{\omega_2}_\gothg + \pair{\LieD{\hat{\omega}_1}{a}}{\abar_2}_V - \int_{\Mbound} [\surfAdStar{\omega_1}(\alpha,\omega_2) \\
&  +(-1)^{c+1}  \surfPhiDTA{a}(\omega_2,\abar_1) - \surfLie{\hat{\omega}_2}{a}{\abar_1} + \surfLie{\hat{\omega}_1}{a}{\abar_2}  ]\bound{},\\
& = \pair{\adSdual{\elS{1}} \elSdual{}}{\elS{2}}_\goths - \int_{\Mbound} \bound{ \surfAdSstar{\elS{1}}(\alpha,a,\omega_2,\abar_2)},
\end{align*}																	
which concludes the proof.

\subsection{Proof of Theorem \ref{theorem:Hamiltonian}}\label{proof:Hamiltonian}
\newcommand{\dOm}{\delta \omega}
\newcommand{\dMu}{\delta \mu}

First we start by the variational derivatives of (\ref{eq:reduced_Lagrangian}) with respect to its variables.
The variational derivative $\delta_{\mu}L_k :=\varD{L_k}{\mu}$ of $L_k$ with respect to $\mu \in V^*$ is the element of $V$ that satisfies for any $\dMu \in V^*$, 
$
\pair{\delta_{\mu} L_k}{\dMu}_{V^*} = \pair{\dMu}{\delta_{\mu} L_k}_{V} =  \depsLine{\epsilon = 0} L_k(\omega_v,\mu+ \epsilon \dMu).
$
By rewriting (\ref{eq:reduced_Lagrangian}) as
$
L_k(\omega_v,\mu) = \int_M \half (*\mu) \wedge \omega_v \wedge * \omega_v = \int_M \half *(\omega_v \wedge *\omega_v) \wedge \mu,
$
then one has that
\begin{equation*}
\delta_{\mu} L_k =  \half *(\omega_v \wedge *\omega_v),
\end{equation*}
by definition of the variational derivative.

The variational derivative $\delta_{\omega_v}L_k :=\varD{L_k}{\omega_v}$ of $L_k$ with respect to $\omega_v \in \gothg$ is the element of $\gothgstar$ that satisfies for any $\dOm \in \gothg$
$
\pair{\delta_{\omega_v} L_k}{\dOm}_\gothg = \depsLine{\epsilon = 0} L_k(\omega_v + \epsilon \dOm,\mu).
$
By rewriting (\ref{eq:reduced_Lagrangian}) as
$
L_k(\omega_v,\mu) = \int_M \half (*\mu) \omega_v \wedge * \omega_v = \int_M \half \mNmOne(*\mu) *\omega_v \wedge  \omega_v ,
$
and the observation that $L_k$ is quadratic in $\omega_v$, it follows, again simply by definition of the variational derivative, that
\begin{equation*}
\delta_{\omega_v} L_k =  \mNmOne (*\mu)*\omega_v.
\end{equation*}

Now define the conjugate momentum variable $\alpha \in \gothgstar$ by $\alpha := \delta_{\omega_v} L_k =  \mNmOne (*\mu)*\omega_v$.
The partial Legendre transform of $L_k$ is by definition the change of variables, and its corresponding inverse,
\begin{equation*}
(\omega_v,\mu) \mapsto (\alpha,\mu) = (\mNmOne (*\mu)*\omega_v,\mu), \qquad \qquad (\alpha,\mu) \mapsto (\omega_v,\mu) =(\frac{*\alpha}{*\mu},\mu).
\end{equation*}
Using (\ref{eq:omega_v_def}), we have that physically $\alpha = (*\mu) \vf = \rho \vf$ is the 1-form corresponding to the momentum of the fluid.

Since the Legendre transformation is a diffeomorphism, then $L_k$ is hyper-regular, and one can define the Hamiltonian $\map{H_k}{\gothgstar\times V^*}{\bb{R}}$ 
\begin{equation}\label{eq:Legendre_Ham}
H_k(\alpha,\mu) := \pair{\alpha}{\omega_v}_\gothg - L_k(\omega_v,\mu),
\end{equation}
which could be expressed as
\begin{align*}
H_k(\alpha,\mu) =& \int_{M} \alpha\wedge\omega_v - \half \mNmOne(*\mu) *\omega_v \wedge  \omega_v,
				= \int_{M} (\alpha - \half \mNmOne(*\mu) *\omega_v) \wedge  \omega_v,\\
				=& \int_{M} (\alpha - \half \alpha) \wedge  \frac{*\alpha}{*\mu},
				= \int_{M} \frac{1}{2(*\mu)} \alpha \wedge  *\alpha.
\end{align*}

The variational derivative of $H$ with respect to $\alpha$ follows immediately by applying the chain rule to the construction in (\ref{eq:Legendre_Ham})
\begin{equation}
\varD{H_k}{\alpha} = \omega_v + \pair{\alpha}{\varD{\omega_v}{\alpha}}_\gothg - \pair{\varD{L_k}{\omega_v}}{\varD{\omega_v}{\alpha}}_\gothg = \omega_v = \frac{*\alpha}{*\mu}.
\end{equation}
Similarly, the variational derivative of $H$ with respect to $\mu$ is given by
\begin{equation*}
\varD{H_k}{\mu} = - \varD{L_k}{\mu} = - \half *(\omega_v \wedge *\omega_v) = -\frac{1}{2(*\mu)^2} *(\alpha\wedge *\alpha) = -\frac{1}{2(*\mu)^2} \iota_{\hat{\alpha}} \alpha,
\end{equation*}
where identity (\ref{eq:identity_Hodge_int_product}) was used. This concludes the proof.


\subsection{Proof of Proposition \ref{proposition:vel_represent_maps}}\label{proof:vel_represent_maps}
For notational simplicity in this proof we denote $*\mu = *\muT$ by $\rho$ when needed.
The derivation of the pushforward map $(\Phi_*)_x$ follows from the rate of change of $\alpha = \rho \vf$
\begin{equation}
\dot\alpha = \dot\rho \vf + \rho \dot\vf \implies  \dot\vf = \frac{\dot\alpha - \dot\rho \vf}{\rho} = \frac{\fAlpha - (*\fMu) \vf}{*\muT},
\end{equation} 
where $\dot \rho$ is the density of the top-form $\dot \mu = \dot \muT$.
This concludes the proof of (\ref{eq:pushfwd_vel_map}).

The pullback map $(\Phi^*)_x$ is defined implicitly by
\begin{equation*}
\pair{(\fAlpha,\fMu)}{(\Phi^*)_x\eXt}_\goths= \pair{(\Phi_*)_x(\fAlpha,\fMu)}{\eXt}_\goths = \pair{(\fVt,\fMuT)}{\eXt}_\goths.
\end{equation*}
Using the duality pairing definition (\ref{eq:duality_s}), we have that
\begin{align*}
\pair{(\fAlpha,\fMu)}{\eX}_\goths =& \int_{M} \fVt\wedge e_{\vf} + \fMuT \wedge e_{\muT} = \int_{M} \left(\frac{\dot\alpha - (*\fMu) \vf}{\rho}\right) \wedge e_{\vf} + \fMu \wedge e_{\muT}\\
													=& \int_{M} \frac{\dot\alpha}{\rho} \wedge e_{\vf} - *\fMu \frac{\vf}{\rho} \wedge e_{\vf}+ \fMu \wedge e_{\muT}
													= \int_{M} \dot\alpha \wedge \frac{e_{\vf}}{\rho}- *\fMu \wedge \frac{\vf\wedge e_{\vf}}{\rho} + \fMu \wedge e_{\muT}\\
													=& \int_{M} \dot\alpha \wedge \frac{e_{\vf}}{\rho}- \fMu \wedge \frac{*(\vf\wedge e_{\vf})}{\rho} + \fMu \wedge e_{\muT}
													= 	\pair{(\fAlpha,\fMu)}{\left(\frac{e_{\vf}}{\rho},e_{\muT} - \frac{*(\vf\wedge e_{\vf})}{\rho} \right)}_\goths,																		
\end{align*}
which concludes the proof of (\ref{eq:pullback_vel_map}).
\newcommand{\Xt}{X_{\vf}}
\newcommand{\XtTilde}{\tilde{X}_{\vf}}

Finally the proof of the Lie-Poisson structure (\ref{eq:JxT_def}) is as follows.
Let $\eX$ be given by (\ref{eq:pullback_vel_map}). Then using (\ref{eq:J_x_2}) one has that
\begin{equation}\label{eq:proof_Jx_1}
\TwoVec{\fAlpha}{\fMu} = J_x\eX = \TwoVec{-\LieD{\Xt}{(\alpha)} - \divr{\Xt}\alpha - \rho \extd (e_{\muT} - \frac{*(\vf\wedge e_{\vf})}{\rho})}{- \extd(\rho \frac{e_{\vf}}{\rho})},
\end{equation}
where $\Xt := (\frac{e_{\vf}}{\rho})^\wedge$ is the vector field corresponding to the $n-1$ form $\frac{e_{\vf}}{\rho} = e_\alpha$ which is defined such that
$\iota_{\Xt} \volF = \frac{e_{\vf}}{\rho} = *\XtTilde.$

The first row in (\ref{eq:proof_Jx_1}) could be massaged as follows
\begin{align}
\fAlpha 	&= -\LieD{\Xt}{(\rho \vf)} - \divr{\Xt}\rho \vf - \rho \extd e_{\muT} +\rho \extd *(\vf\wedge \frac{e_{\vf}}{\rho})\nonumber\\
				&= -\rho\LieD{\Xt}{\vf} -\LieD{\Xt}{(\rho)}\vf - \divr{\Xt}\rho \vf - \rho \extd e_{\muT} +\rho \extd *(\vf\wedge *\XtTilde)\nonumber\\
				&= -\rho\extd \iota_{\Xt}{\vf} -\rho\iota_{\Xt}{\extd\vf} -(\LieD{\Xt}{\rho}+ \divr{\Xt}\rho) \vf - \rho \extd e_{\muT} +\rho \extd *(\XtTilde\wedge *\vf)\nonumber\\
				&= -\rho\iota_{\Xt}{\extd\vf} -(\LieD{\Xt}{\rho}+ \divr{\Xt}\rho) \vf - \rho \extd e_{\muT}\label{eq:alphaDot_proof},									
\end{align}
where the Leibniz rule for the Lie derivative and Cartan's homotopy formula were used, and the last equality follows from $-\iota_{\Xt}{\vf}+*(\XtTilde\wedge *\vf) = - \iota_{\Xt}{\vf} + \iota_{\Xt}{\vf} =0$, using identity (\ref{eq:identity_Hodge_int_product}).
The second row in (\ref{eq:proof_Jx_1}), and consequently the second row in (\ref{eq:Jx_tilde}), can be written simply as
\begin{equation}\label{eq:muDot_proof}
\fMu = - \extd e_{\vf} = \fMuT.
\end{equation}

By substituting the expressions (\ref{eq:alphaDot_proof}) in the expression of $\fVt$ in (\ref{eq:pushfwd_vel_map}), we have that
\begin{equation*}
\fVt = \frac{1}{\rho}(\fAlpha - \dot{\rho}\vf) = -\iota_{\Xt}{\extd\vf} -\frac{1}{\rho}\underbrace{(\dot{\rho}+\LieD{\Xt}{\rho}+ \divr{\Xt}\rho)}_\text{=0} \vf - \extd e_{\muT} =  -\iota_{\Xt}{\extd\vf} - \extd e_{\muT},
\end{equation*}
where the term in the parenthesis vanishes as a consequence of the mass continuity equation (\ref{eq:continuity_density}).
Thus, finally using 
$\iota_{\Xt}{\extd\vf}= \iota_{(\frac{e_{\vf}}{\rho})^\wedge}{\extd\vf} = \frac{1}{\rho}\iota_{\hat{e}_{\vf}}{\extd\vf},$
the proof of (\ref{eq:Jx_tilde}) is concluded.

\subsection{Proof of Proposition \ref{prop:kinetic_Ham_new}}\label{{proof:kinetic_Ham_new}}

The new functional $\tilde{H}_k\stateV$ is defined by $\tilde{H}_k = H_k \circ\Phi^{-1}$. By substituting $\Phi^{-1}\stateV = ((*\muT)\vf,\mu)$ in (\ref{eq:reduced_Hamiltonian}) as
$$\tilde{H}_k\stateV = H_k (\Phi^{-1}\stateV) =  \int_M \frac{1}{2(*\mu)} (*\mu)\vf \wedge * (*\mu)\vf = \int_M \half (*\mu) \vf \wedge * \vf.$$
For notational simplicity, we drop the \textit{tilde} from the mass form $\mu = \muT\in V^*$ as well as the tilde over the new energy functional and denote it by $H_k\stateV$.

The variational derivative of $H_k$ with respect to $\vf\in \gothgstar$ can be found out by rewriting (\ref{eq:new_Hamiltonian}) as 
$H_k\stateV = \int_M \half \vf \wedge ((*\mu)* \vf),$
where we have that $\delta_{\vf} H_k  = (*\mu) * \vf$ from the fact that $H_k$ is quadratic in $\vf$.
Furthermore, we can express $\delta_{\vf} H_k$ as 
$$\delta_{\vf} H_k = (*\mu) * \vf = (*\mu) \iota_v\volF = \iota_v (*\mu \volF) = \iota_v \mu.$$

Similarly, the variational derivative of $H_k$ with respect to $\mu \in V^*$ can be found out by rewriting (\ref{eq:new_Hamiltonian}) as 
$H_k(\vf,\mu) = \int_M \half (*\mu) \wedge (\vf \wedge * \vf) = \int_M \mu \wedge \half*(\vf \wedge * \vf),$
where $$ \delta_{\mu} H_k = \half*(\vf \wedge * \vf)$$ follows by definition of variational derivative.
Using identity (\ref{eq:identity_Hodge_int_product}) allows one to write $ \delta_{\mu} H_k$ as in (\ref{eq:varDeriv_stateV}), which concludes the proof.

\subsection{Proof of Theorem \ref{theorem:new_power_balance_H_k}} \label{proof:new_power_balance_H_k}

In this proof we denote by $\eX = (\eAlpha,\eMu)$ and by $\eXt = (\eVt,\eMuT)$.
By substituting $\eX$ defined by (\ref{eq:pullback_vel_map}) and $\state = (*\mu \vf,\mu)$ in the energy balance expression in terms of $\state$ given by (\ref{eq:energy_balance_old})
\begin{align*}
\dot{H}_k =& \int_{\Mbound} { \bound{\left(\frac{\iota_{X_\alpha}(\alpha)}{*\mu} + e_\mu \right)} \wedge \bound{-(*\mu)e_\alpha} }  + \int_{M} \fstr \wedge e_\alpha,\\
			=& \int_{\Mbound} \bound{\left(\iota_{X_\alpha}(\vf) + e_{\muT} - \frac{*(\vf\wedge e_{\vf})}{*\mu}\right)} \wedge - \bound{e_\vf}  + \int_{M} \fstr\wedge \frac{e_\vf}{*\mu} ,
\end{align*}
where $X_\alpha$ is defined such that $\iota_{X_\alpha}\volF = e_\alpha = \frac{e_\vf}{*\mu}$. Consequently, we have that $\frac{e_\vf}{*\mu} = *\tilde{X}_\alpha$.
The term $\frac{*(\vf\wedge e_{\vf})}{*\mu}$ can be rewritten using identity (\ref{eq:identity_Hodge_int_product}) as
$\frac{*(\vf\wedge e_{\vf})}{*\mu} = *(\vf\wedge \frac{e_{\vf}}{*\mu}) = *(\vf\wedge *\tilde{X}_\alpha) =  *(\tilde{X}_\alpha \wedge *\vf) = \iota_{X_\alpha} \vf.$
Therefore, we have that
\begin{align*}
\dot{H}_k =& \int_{\Mbound} \bound{\left(\iota_{X_\alpha}(\vf) + e_{\muT} - \frac{*(\vf\wedge e_{\vf})}{*\mu}\right)} \wedge - \bound{e_\vf}  + \int_{M} \fstr\wedge \frac{e_\vf}{*\mu},\\
			=& \int_{\Mbound} { \bound{e_{\muT}} \wedge - e_\vf\bound{}} + \int_{M} \fstr\wedge \frac{e_\vf}{*\mu}.
\end{align*}
Therefore, $e_\vf = \eVt$ and $e_{\muT} = \eMuT$ are given in (\ref{eq:varDeriv_stateV}) (where $\muT =\mu$ is used).
Moreover, the term $\frac{e_\vf}{*\mu}$ could be rewritten using (\ref{eq:varDeriv_stateV}) as
$\frac{e_\vf}{*\mu} = \frac{\eVt}{*\mu} = \frac{1}{*\mu} \iota_v\mu = \iota_v(\frac{\mu}{*\mu}) = \iota_v\volF,$
which concludes the proof.

\bibliographystyle{abbrv}
\bibliography{Paper_References}

\end{document}